ANNALES
DE L'INSTITUT
HENRI
POINCARÉ
PROBABILITÉS
ET STATISTIQUES

# On homogenization of space-time dependent and degenerate random flows II

## Rémi Rhodes

*Laboratoire d'Analyse Topologie Probabilités, Université de Provence, 39 rue Joliot-Curie, 13453 Marseille Cedex 13, France. E-mail: rhodes@cmi.univ-mrs.fr*

**Abstract.** We study the long time behavior (homogenization) of a diffusion in random medium with time and space dependent coefficients. The diffusion coefficient may degenerate. In *Stochastic Process. Appl.* (2007) (to appear), an invariance principle is proved for the critical rescaling of the diffusion. Here, we generalize this approach to diffusions whose space-time scaling differs from the critical one.

**Résumé.** Nous étudions le comportement asymptotique (homogénéisation) d'une diffusion en milieu aléatoire avec des coefficients dépendant du temps et de l'espace, pour laquelle le coefficient de diffusion peut dégénérer. Dans *Stochastic Process. Appl.* (2007) (to appear), un principe d'invariance est établi pour le changement d'échelle critique de la diffusion. Ici, une généralisation de cette approche est proposée pour différents changements d'échelle possibles.

## 1. Introduction

We aim at studying the long time behavior, as $\varepsilon \to 0$, of the diffusion process with random coefficients (the parameter $\omega$ below stands for this randomness) defined by

$$X_t^{\varepsilon,\omega} = x + \varepsilon^{-\beta} \int_s^t b\left(\frac{r}{\varepsilon^\alpha}, \frac{X_r^{\varepsilon,\omega}}{\varepsilon^\beta}, \omega\right) \mathrm{d}r + \int_s^t \sigma\left(\frac{r}{\varepsilon^\alpha}, \frac{X_r^{\varepsilon,\omega}}{\varepsilon^\beta}, \omega\right) \mathrm{d}B_r \tag{1}$$

and then at identifying the limit of the solutions of the random parabolic differential equations (PDE)

$$\partial_t z_{\varepsilon,\omega}(x,t) = \frac{1}{2} \text{trace}\left[a\left(\frac{t}{\varepsilon^\alpha}, \frac{x}{\varepsilon^\beta}, \omega\right) \Delta_{xx} z_{\varepsilon,\omega}(x,t)\right] + \varepsilon^{-\beta} b\left(\frac{t}{\varepsilon^\alpha}, \frac{x}{\varepsilon^\beta}, \omega\right) \cdot \nabla_x z_{\varepsilon,\omega}(x,t)$$

$$+ [\varepsilon^{-\beta} c + d]\left(\frac{t}{\varepsilon^\alpha}, \frac{x}{\varepsilon^\beta}, \omega\right) z_{\varepsilon,\omega}(x,t), \quad (x,t) \in \mathbb{R}^d \times \mathbb{R}_+ \tag{2}$$

with initial condition $z_{\varepsilon,\omega}(x,0) = f(x)$. $\alpha, \beta$ are two strictly positive parameters. The coefficients $a, b, \sigma, c, d$ are stationary ergodic random fields with respect to space and time variables. We shall see that there are different possible limits depending on $\alpha = 2\beta, \alpha < 2\beta$ or $\alpha > 2\beta$. More precisely, we suppose that $a = \sigma\sigma^*$ and the generator of the diffusion process could be rewritten in divergence form as

$$L^{\varepsilon,\omega} = \left(\frac{1}{2}\right)e^{2V(x/\varepsilon^\beta,\omega)}\mathrm{div}_x\left(e^{-2V(x/\varepsilon^\beta,\omega)}[a+H]\left(\frac{t}{\varepsilon^\alpha},\frac{x}{\varepsilon^\beta},\omega\right)\nabla_x\right). \tag{3}$$

Here $V$ and $H$ are also stationary random fields and $H$ is antisymmetric. Assumptions will be stated rigorously in the next section. We will prove that, in probability with respect to $\omega$,

$$\lim_{\varepsilon\to 0} z_{\varepsilon,\omega}(x,t) = \overline{z}(t,x), \tag{4}$$

where $\overline{z}$ is the solution of a deterministic equation of the type

$$\partial_t \overline{z}(x,t) = \mathrm{trace}[A\Delta_{xx}\overline{z}](x,t) + C\cdot\nabla_x\overline{z}(x,t) + U\overline{z}(x,t) \tag{5}$$

with initial condition $\overline{z}(x,0) = f(x)$. $A, C, U$ are deterministic coefficients, the so-called effective coefficients, and depend on the considered case. As in the previous paper [18], we do not assume any non-degeneracy of the diffusion matrix $a$ (several examples are introduced in [18]).

Roughly speaking, the case $\alpha = 2\beta$ corresponds to the critical scaling for the diffusion and space and time variables have to be homogenized simultaneously. This job was carried out in [18]. In case $\alpha < 2\beta$, the space variables are moving faster than the time variable so that the homogenization procedure has to be performed first in space and then in time, and vice versa in the case $\alpha > 2\beta$.

Both this paper and [18] follow a series of works on this topic. Let us sum up the methodological approach of this issue. Briefly, the time dependence of the process brings about a strong time degeneracy of the underlying Dirichlet form, which satisfies no sector condition, even weak. To face such a degeneracy, analytical methods, more precisely compactness methods, carried through the homogenization procedure for a uniformly elliptic diffusion matrix $a$ in periodic media [1, 20] or in random media [4]. In [8], critical scaling is considered ($\alpha = 2$ and $\beta = 1$) and an annealed invariance principle is proved for the diffusion (1) by means of probabilistic tools under the assumptions $V = 0$ and $a = \mathrm{Id}$. This result is extended in [10] to the case when the stream matrix $H$ satisfies a certain integrability condition of the spatial energy spectrum instead of boundedness of the coefficients. A quenched version of the invariance principle is stated in [5] for an arbitrary time and space dependent diffusion coefficient provided that it satisfies a strong uniform non-degeneracy assumption. The particularity of these works lies in their intensive use of regularity properties of the heat kernel to deal with the time degeneracy of the Dirichlet form.

In this paper, we additionally consider possibly degenerate diffusion coefficients $a$. This prevents us from using both compactness methods (the Poincaré inequality is lacking) and regularity properties of the heat kernel (no uniform ellipticity or even hypo-ellipticity of the matrix $a$ is assumed). The main tools of the proof are the control of the matrix $a$ by a time independent one $\tilde{a}$ and a commutativity argument of the unbounded operators associated to $\tilde{a}$ and to the time evolution.

This article is structured as follows. In Section 2, we introduce notations and present our results. In Section 3, we set out the main properties of the involved stochastic processes. Section 4 is devoted to constructing the correctors. We explain how to get round the lack of regularity of the correctors in Section 5 to establish the Itô formula. Then in Section 6, an ergodic theorem is proved, which allows us to carry through the homogenization procedure (Section 7). The tightness of $X^{\varepsilon,\omega}$ is studied in Section 8 by means of the Garsia–Rodemich–Rumsey inequality. The limit PDE 5 is identified in Section 9 with the help of the Girsanov transform.

## 2. Notations, setup and main results

The setup remains the same as in [18]. It is reminded for the reader's convenience.

**Definition 2.1.** Let $(\Omega,\mathcal{G},\mu)$ be a probability space and $\{\tau_{t,x}; (t,x)\in\mathbb{R}\times\mathbb{R}^d\}$ a stochastically continuous group of measure preserving transformations acting ergodically on $\Omega$:

(1)  $\forall A \in \mathcal{G}, \forall (t,x) \in \mathbb{R} \times \mathbb{R}^d, \, \mu(\tau_{t,x}A) = \mu(A),$

(2)  *If for any* $(t,x) \in \mathbb{R} \times \mathbb{R}^d, \, \tau_{t,x}A = A$ *then* $\mu(A) = 0 \, or \, 1,$

(3)  *For any measurable function* $\mathbf{g}$ *on* $(\Omega, \mathcal{G}, \mu)$, *the function* $(t,x,\omega) \mapsto \mathbf{g}(\tau_{t,x}\omega)$ *is measurable on* $(\mathbb{R} \times \mathbb{R}^d \times \Omega, \mathcal{B}(\mathbb{R} \times \mathbb{R}^d) \otimes \mathcal{G})$.

In what follows we will use the bold type to denote a function $\mathbf{g}$ from $\Omega$ into $\mathbb{R}$ (or more generally into $\mathbb{R}^n$, $n \geq 1$) and the unbold type $g(t,x,\omega)$ to denote the associated representation mapping $(t,x,\omega) \mapsto \mathbf{g}(\tau_{t,x}\omega)$. The space of square integrable functions on $(\Omega, \mathcal{G}, \mu)$ is denoted by $L^2(\Omega)$, the usual norm by $|\cdot|_2$ and the corresponding inner product by $(\cdot, \cdot)_2$. Then, the operators on $L^2(\Omega)$ defined by $T_{t,x}\mathbf{g}(\omega) = \mathbf{g}(\tau_{t,x}\omega)$ form a strongly continuous group of unitary maps in $L^2(\Omega)$. Each function $\mathbf{g}$ in $L^2(\Omega)$ defines in this way a stationary ergodic random field on $\mathbb{R}^{d+1}$. The group possesses $d+1$ generators defined for $i = 1, \ldots, d$, by

$$D_i\mathbf{f} = \frac{\partial}{\partial x_i}T_{0,x}\mathbf{f}|_{(t,x)=0} \quad \text{and} \quad D_t\mathbf{f} = \frac{\partial}{\partial t}T_{t,0}\mathbf{f}|_{(t,x)=0},$$

which are closed and densely defined. Denote by $\mathcal{C}$ the dense subset of $L^2(\Omega)$ defined by

$$\mathcal{C} = \mathrm{Span}\{\mathbf{f} * \varphi; \mathbf{f} \in L^2(\Omega), \varphi \in C_c^\infty(\mathbb{R}^{d+1})\}, \quad \text{with } \mathbf{f} * \varphi(\omega) = \int_{\mathbb{R}^{d+1}} \mathbf{f}(\tau_{t,x}\omega)\varphi(t,x)\,dt\,dx,$$

where $C_c^\infty(\mathbb{R}^{d+1})$ is the set of smooth functions on $\mathbb{R}^{d+1}$ with a compact support. Remark that $\mathcal{C} \subset \mathrm{Dom}(D_i)$ and $D_i(\mathbf{f} * \varphi) = -\mathbf{f} * \frac{\partial \varphi}{\partial x_i}$. This last quantity is also equal to $D_i\mathbf{f} * \varphi$ if $\mathbf{f} \in \mathrm{Dom}(D_i)$.

Consider now measurable functions $\boldsymbol{\sigma}, \widetilde{\boldsymbol{\sigma}}, \mathbf{H} : \Omega \to \mathbb{R}^{d \times d}$ and $\mathbf{V}, \mathbf{c}, \mathbf{d} : \Omega \to \mathbb{R}$. Assume that $\mathbf{H}$ is antisymmetric. Define $\mathbf{a} = \boldsymbol{\sigma}\boldsymbol{\sigma}^*$ and $\widetilde{\mathbf{a}} = \widetilde{\boldsymbol{\sigma}}\widetilde{\boldsymbol{\sigma}}^*$. The function $\mathbf{V}$ does not depend on time, that means $\forall t \in \mathbb{R}$, $T_{t,0}\mathbf{V} = \mathbf{V}$.

### Assumption 2.2 (Regularity of the coefficients).

- *Assume that* $\forall i,j,k,l = 1, \ldots, d, \; \mathbf{a}_{ij}, \widetilde{\mathbf{a}}_{ij}, \mathbf{V}, \mathbf{H}_{ij}, D_l\mathbf{a}_{ij}, \, and \, D_l\widetilde{\mathbf{a}}_{ij} \in \mathrm{Dom}(D_k)$.
- $\mathbf{c} = e^{2\mathbf{V}}\sum_{i,j} D_i(e^{-2\mathbf{V}}\widetilde{\boldsymbol{\sigma}}_{ij}\mathbf{f}_j)$ *for some function* $\mathbf{f} \in L^\infty(\Omega; \mathbb{R}^d)$.
- *Define, for* $i = 1, \ldots, d,$

$$
\begin{aligned}
\mathbf{b}_i(\omega) &= \sum_{j=1}^d \left( \frac{1}{2}D_j\mathbf{a}_{ij}(\omega) - \mathbf{a}_{ij}D_j\mathbf{V}(\omega) + \frac{1}{2}D_j\mathbf{H}_{ij}(\omega) \right), \\
\widetilde{\mathbf{b}}_i(\omega) &= \sum_{j=1}^d \left( \frac{1}{2}D_j\widetilde{\mathbf{a}}_{ij}(\omega) - \widetilde{\mathbf{a}}_{ij}D_j\mathbf{V}(\omega) \right)
\end{aligned}
\tag{6}
$$

*and assume that the applications* $(t,x) \mapsto b_i(t,x,\omega)$, $(t,x) \mapsto c(t,x,\omega)$, $(t,x) \mapsto \sigma(t,x,\omega)$, $x \mapsto DV(x,\omega)$ *are globally Lipschitz (for each fixed* $\omega \in \Omega$) *and the application* $(t,x) \mapsto d(t,x,\omega)$ *is continuous. Moreover, the coefficients* $\boldsymbol{\sigma}, \mathbf{a}, \mathbf{b}, \widetilde{\boldsymbol{\sigma}}, \mathbf{V}, \mathbf{H}, \mathbf{c}, \mathbf{d}$ *are bounded by a constant* $K$. *(In particular, this ensures existence and uniqueness of a global solution of SDE (1).)*

Here is the main assumption of this paper:

### Assumption 2.3 (Control of the coefficients).

- $\widetilde{\boldsymbol{\sigma}}$ *does not depend on time (i.e.* $\forall t \in \mathbb{R}$, $T_t\widetilde{\boldsymbol{\sigma}} = \widetilde{\boldsymbol{\sigma}}$). *As a consequence, the matrix* $\widetilde{\mathbf{a}}$ *does not depend on time either.*
- $\mathbf{H}, \mathbf{a}, \mathbf{c} \in \mathrm{Dom}(D_t)$ *and there exist five positive constants* $m(>0), M, C_1^H, C_2^H, C_2^a$ *such that,* $\mu$ *a.s.,*

$$m\widetilde{\mathbf{a}} \leq \mathbf{a} \leq M\widetilde{\mathbf{a}}, \tag{7}$$

$$|\mathbf{H}| \leq C_1^H\widetilde{\mathbf{a}}, \qquad |D_t\mathbf{H}| \leq C_2^H\widetilde{\mathbf{a}} \quad and \quad |D_t\mathbf{a}| \leq C_2^a\widetilde{\mathbf{a}}, \tag{8}$$

where $|\mathbf{A}|$ stands for the symmetric positive square root of $\mathbf{A}$, i.e. $|\mathbf{A}| = \sqrt{\mathbf{A}\mathbf{A}^*}$.

For instance, if the matrix $\mathbf{a}$ is uniformly elliptic and bounded, $\widetilde{\boldsymbol{\sigma}}$ can be chosen as equal to the identity matrix $\mathbf{Id}$ and then (8) $\Leftrightarrow \mathbf{H}, D_t\mathbf{H}$ and $D_t\mathbf{a} \in L^\infty(\Omega)$.

Let us now set out the ergodic properties of this framework:

**Assumption 2.4 (Ergodicity).** *Let us consider the operator $\widetilde{\mathbf{S}} = (1/2)\mathrm{e}^{2\mathbf{V}}\sum_{i,j=1}^d D_i(\mathrm{e}^{-2\mathbf{V}}\widetilde{\mathbf{a}}_{ij}D_j)$ with domain $\mathcal{C}$. From Assumption 2.2, we can consider its Friedrich extension (see [6], Chapter 3, Section 3), which is still denoted $\widetilde{\mathbf{S}}$. Assume that each function $\mathbf{f} \in \mathrm{Dom}(\widetilde{\mathbf{S}})$ satisfying $\widetilde{\mathbf{S}}\mathbf{f} = 0$ must be $\mu$ almost surely equal to some function that is invariant under space translations.*

The reader is referred to [18] for examples in which Assumptions 2.2, 2.3 and 2.4 are satisfied.

Even if it means adding to $\mathbf{V}$ a constant (and this does not change the drift $\mathbf{b}$, see (6)), we assume that $\int \mathrm{e}^{-2\mathbf{V}}\,\mathrm{d}\mu = 1$. Thus we can define a new probability measure on $\Omega$ by

$$\mathrm{d}\pi(\omega) = \mathrm{e}^{-2\mathbf{V}(\omega)}\,\mathrm{d}\mu(\omega).$$

We now consider a standard $d$-dimensional Brownian motion defined on a probability space $(\Omega', \mathcal{F}, \mathbb{P})$ and the diffusion in random medium given by (1). Under Assumptions 2.2, 2.3 and 2.4, we claim:

**Theorem 2.5.** *The law of the process $X^{\varepsilon,\omega}$ converges weakly in $\mu$-probability in $C([0,T];\mathbb{R}^d)$ to the law of a Brownian motion with a certain covariance matrix (see (50)).*

This result puts us in position to describe the long time behavior of PDE (2):

**Theorem 2.6.** *For any bounded continuous function $f:\mathbb{R}^d \to \mathbb{R}$, the solution $z_{\varepsilon,\omega}(x,t)$ of PDE (2) with initial condition $z_{\varepsilon,\omega}(x,0) = f(x)$ converges in $\mu$-probability towards $\bar{z}(x,t)$, where $\bar{z}$ is the unique viscosity solution of the deterministic PDE (5) with initial condition $\bar{z}(x,0) = f(x)$ (see (49) and (50) for a description of $A, C, U$). More precisely, for any $(x,t) \in \mathbb{R}^d \times \mathbb{R}_+$ and $\delta > 0$*

$$\mu(\{\omega; |z_{\varepsilon,\omega}(x,t) - \bar{z}(x,t)| \geq \delta\}) \to 0$$

*as $\varepsilon$ tends to 0. Furthermore, the coefficients $A, C, U$ only depend on the case $\alpha - 2\beta < 0$, $\alpha - 2\beta > 0$ and $\alpha - 2\beta = 0$.*

The reader can refer to [2] for a description of the viscosity solution theory of PDEs.

**Remark 2.7.** *Let us be more explicit about the last statement of Theorem 2.6. Because of the possible degeneracies of the diffusion matrix $\mathbf{a}$, the coefficients $A, C, U$ are defined by limits involving solutions of parameterized PDEs (see Section 9). However, the uniformly elliptic setting provides us with exact problems (or PDEs) to compute $A, C, U$. The reader is referred to [1], Chapter 2, Section 1, or [16, 20] for the periodic case and to [4] for the random case.*

To get an idea of the scaling effect, let us focus on the 1-dimensional periodic case. In this setting, the matrix $A$ is given by (see [20], Proposition 2.2)

$$A = \int_{t\in[0,1]}\int_{x\in[0,1]} \mathbf{a}(t,x)(D\mathbf{v}(t,x)+1)\,\mathrm{d}x\,\mathrm{d}t, \tag{9}$$

where the function $\mathbf{v}$ solves one of the following problems, depending on the values of the parameters $\alpha$ and $\beta$.

• Case $\alpha - 2\beta < 0$: $\mathbf{v}$ is the unique solution of the parameter dependent elliptic problem:

$$\begin{cases} -D_x(\mathbf{a}(t,x)(D_x\mathbf{v}(t,x)+1)) = 0, \\ \mathbf{v}(t,\cdot) \text{ is 1-periodic}, \quad \int_{x\in[0,1]} \mathbf{v}(t,x)\,\mathrm{d}x = 0, \quad t \in [0,1]. \end{cases}$$

It is readily seen that $D_x\mathbf{v}(t,x) = C(t)/\mathbf{a}(t,x) - 1$ where $C$ is a function of $t$. Since $\mathbf{v}$ is $x$-periodic, $\int_{x\in[0,1]} D_x\mathbf{v}(t,x)\,dx = 0$. Consequently, it comes out that $C(t) = (\int_{x\in[0,1]} \frac{1}{\mathbf{a}(t,x)}\,dx)^{-1}$ so that, from *(9)*, we obtain explicitly the effective diffusivity

$$A = \int_{t\in[0,1]} \left(\int_{x\in[0,1]} \frac{1}{\mathbf{a}(t,x)}\,dx\right)^{-1} dt.$$

- Case $\alpha - 2\beta > 0$: $\mathbf{v}$ is the unique solution of the elliptic problem:

$$\begin{cases} -D_x(\breve{\mathbf{a}}(x)(D_x\mathbf{v}(x)+1)) = 0, \\ \mathbf{v} \text{ is 1-periodic,} \quad \int_{x\in[0,1]} \mathbf{v}(x)\,dx = 0, \end{cases}$$

where $\breve{\mathbf{a}}(x) = \int_{t\in[0,1]} \mathbf{a}(t,x)\,dt$. As in the previous case, solving this problem raises no particular difficulty and we get $D_x\mathbf{v}(x) = C/\breve{\mathbf{a}}(x) - 1$ where $C$ is a constant, which exactly matches $(\int_{x\in[0,1]} \frac{1}{\breve{\mathbf{a}}(x)}\,dx)^{-1}$ thanks to the periodicity of $\mathbf{v}$ again. Then formula *(9)* reads

$$A = \left(\int_{x\in[0,1]} \frac{1}{\int_{t\in[0,1]} \mathbf{a}(t,x)\,dt}\,dx\right)^{-1}.$$

*Note that, in comparison with the previous case, the role of time and space variables have been, in a way, inverted.*

- Case $\alpha - 2\beta = 0$: $\mathbf{v}$ is the unique solution of the parabolic problem:

$$\begin{cases} D_t\mathbf{v}(t,x) - D_x(\breve{\mathbf{a}}(t,x)(D_x\mathbf{v}(t,x)+1)) = 0, \\ \mathbf{v} \text{ is periodic in time and space,} \quad \int_{t,x\in[0,1]} \mathbf{v}(t,x)\,dt\,dx = 0. \end{cases}$$

*Unfortunately, computations in the critical case $\alpha = 2\beta$ are trickier and do not enable us to give an explicit formula for the effective diffusivity.*

From now on, in order to avoid heavy notations, the superscript $\omega$ of $X^{\varepsilon,\omega}$ is omitted when there is no possible confusion. So we write $X^\varepsilon$ instead of $X^{\varepsilon,\omega}$. Moreover, except where otherwise stated, the process $X^\varepsilon$ starts at time $s = 0$ (see *(1)*).

## 3. Environment as seen from the particle

Let us denote by $\overline{X}^\varepsilon$ the process driven by the following Itô equation

$$\overline{X}_t^\varepsilon = \int_0^t b(\varepsilon^{2\beta-\alpha}r, \overline{X}_r^\varepsilon, \omega)\,dr + \int_0^t \sigma(\varepsilon^{2\beta-\alpha}r, \overline{X}_r^\varepsilon, \omega)\,dB_r. \tag{10}$$

From law uniqueness for solutions of Itô equations with Lipschitz coefficients, the processes $\{\varepsilon^\beta \overline{X}_{t/\varepsilon^{2\beta}}^\varepsilon; t \geq 0\}$ and $\{X_t^\varepsilon; t \geq 0\}$, which both start from 0, have the same law.

We now focus on the *environment as seen from the particle*: this is a process on the probability space $\Omega$ defined for each $\varepsilon > 0$ by

$$Y_t^\varepsilon(\omega) = \tau_{\varepsilon^{2\beta-\alpha}t, \overline{X}_t^\varepsilon}\omega, \tag{11}$$

where the process $\overline{X}^\varepsilon$ starts from $0 \in \mathbb{R}^d$. This section is devoted to proving that it defines a $\Omega$-valued Markov process, which admits $\pi$ as (not necessarily unique) invariant measure. As a consequence, the associated semigroup (see [3] for a thorough study on semigroups) can be extended to a contraction semigroup on $L^p(\Omega)$, for $1 < p < \infty$. Let us additionally mention that, with the help of the Itô formula, the generator of this process is easily identified on $\mathcal{C}$. It coincides with the operator $\mathbf{L} + \varepsilon^{2\beta - \alpha} D_t$, where $\mathbf{L}$ is defined on $\mathcal{C}$ by

$$\mathbf{L} = \frac{e^{2\mathbf{V}}}{2} \sum_{i,j=1}^{d} D_i(e^{-2\mathbf{V}}[\mathbf{a} + \mathbf{H}]_{ij} D_j). \tag{12}$$

Note that $\mathcal{C}$ need not be a core for this operator because of the possible degeneracies of the matrix $\mathbf{a}$.

All these statements are well known in the case of a uniform elliptic diffusion coefficient $\mathbf{a}$ for time independent (see [9, 13, 14, 15]) or time dependent (see [8, 10] and references therein) coefficients. In what follows, we will see that the degenerate case boils down to the uniform elliptic case by means of vanishing viscosity methods.

Let us consider a $(d+1)$-dimensional Brownian motion $B'$ independent of $B$. Up to the end of this section, the couple $(a, X), a \in \mathbb{R}, X \in \mathbb{R}^d$ stands for a $(d+1)$-dimensional vector and the variables $u, v$ denote $(d+1)$-dimensional vectors. Let us define the $(d+1)$-dimensional process $X^{\varepsilon,n}$ as the solution of the following Itô equation

$$dX_t^{\varepsilon,n} = (\varepsilon^{2\beta-\alpha}, b - n^{-1} DV(X_t^{\varepsilon,n}, \omega)) dt + (0, \sigma(X_t^{\varepsilon,n}, \omega) dB_t) + n^{-1/2} dB_t'. \tag{13}$$

From Lemma 3.1 below, the $\mathbb{R}^{d+1}$-valued process $X^{\varepsilon,n}$ admits transition densities $p_{\varepsilon,n}(t, u, v)$ with respect to the Lebesgue measure on $\mathbb{R}^{d+1}$. The transition densities satisfy for any $t > 0$ and $u, v \in \mathbb{R}^{d+1}$

$$\frac{1}{Ct^{(d+1)/2}} \exp\left(-\frac{C|v - u - (\varepsilon^{2\beta-\alpha}t, 0)|^2}{t}\right) \leq p_{\varepsilon,n}(t, u, v) \leq \frac{C}{t^{(d+1)/2}} \exp\left(-\frac{|v - u - (\varepsilon^{2\beta-\alpha}t, 0)|^2}{Ct}\right) \tag{14}$$

for some constant $C > 0$ that only depends on $K, d, n$. In particular, we can now easily adapt the proof in [9], Section 2.3, and prove that the $\Omega$-valued process $Y_t^{\varepsilon,n}(\omega) = \tau_{X_t^{\varepsilon,n}}\omega$, $X^{\varepsilon,n}$ starting from 0 at time $t = 0$, is a Markov process whose generator coincides on $\mathcal{C}$ with

$$\mathbf{L}^{n,\varepsilon} = \frac{e^{2\mathbf{V}}}{2} \sum_{i,j=1}^{d} D_i(e^{-2\mathbf{V}}[\mathbf{a} + \mathbf{H} + n^{-1}\mathbf{Id}]_{ij} D_j) + (2n)^{-1} D_t^2 + \varepsilon^{2\beta-\alpha} D_t, \tag{15}$$

and admits $\pi$ as invariant measure.

It now remains to let the parameter $n$ go to $\infty$. For this purpose, define the subspace $\mathrm{C}(\Omega)$ of $L^\infty(\Omega)$ as follows: $\mathbf{f} \in L^\infty(\Omega)$ belongs to $\mathrm{C}(\Omega)$ if and only if it is bounded and for $\mu$ almost every $\omega \in \Omega$, the function $(t, x) \in \mathbb{R} \times \mathbb{R}^d \mapsto \mathbf{f}(\tau_{t,x}\omega)$ is continuous on $\mathbb{R} \times \mathbb{R}^d$. Let $\bar{\mathrm{C}}(\Omega)$ denote the closure of $\mathrm{C}(\Omega)$ in $L^\infty(\Omega)$ with respect to the usual $L^\infty(\Omega)$-norm. For each function $\mathbf{f} \in \mathrm{C}(\Omega)$, let us then define

$$P_t(\mathbf{f})(\omega) = \mathbb{E}_0[\mathbf{f}(\tau_{\varepsilon^{2\beta-\alpha}t, \overline{X}_t^\varepsilon}\omega)] = \mathbb{E}[\mathbf{f}(Y_t^\varepsilon(\omega))].$$

Obviously, $P_t$ is continuous with respect to the $L^\infty(\Omega)$-norm. Lemma 3.2 below ensures that $P_t(\mathrm{C}(\Omega)) \subset \mathrm{C}(\Omega)$. Then, from the Markov property of the $\mathbb{R}^{d+1}$-valued process $(\varepsilon^{2\beta-\alpha}t, \overline{X}_t^\varepsilon)$ and (18), it is readily seen that $P_s(P_t\mathbf{f}) = P_{s+t}(\mathbf{f})$ so that the family $(P_t)_{t \in \mathbb{R}_+}$ defines a semigroup on $\bar{\mathrm{C}}(\Omega)$.

Let us now consider a function $\mathbf{f} \in \mathrm{C}(\Omega)$ and fix $\omega \in \Omega$. From the Lipschitz property of the coefficients and the continuity of $(t, x) \in \mathbb{R} \times \mathbb{R}^d \mapsto \mathbf{f}(\tau_{t,x}\omega)$, classical arguments of SDE theory imply that $\mathbb{E}_0[\mathbf{f}(\tau_{X_t^{\varepsilon,n}}\omega)]$ converges to $\mathbb{E}_0[\mathbf{f}(\tau_{\varepsilon^{2\beta-\alpha}t, \overline{X}_t^\varepsilon}\omega)]$ as $n$ goes to $\infty$. Together with the boundedness of $\mathbf{f}$, this ensures that $\pi(\mathbb{E}_0[\mathbf{f}(\tau_{X_t^{\varepsilon,n}}\omega)])$ converges to $\pi(\mathbb{E}_0[\mathbf{f}(\tau_{\varepsilon^{2\beta-\alpha}t, \overline{X}_t^\varepsilon}\omega)])$ as $n$ tends to $\infty$. As $\pi(\mathbb{E}_0[\mathbf{f}(\tau_{X_t^{\varepsilon,n}}\omega)]) = \pi(\mathbf{f})$, it comes out $\pi(\mathbb{E}_0[\mathbf{f}(\tau_{\varepsilon^{2\beta-\alpha}t, \overline{X}_t^\varepsilon}\omega)]) = \pi(\mathbf{f})$. Hence $\pi$ is invariant for the semigroup $(P_t)_{t \in \mathbb{R}_+}$, which now clearly extends to a contraction semigroup on $L^p(\Omega, \pi)$ for $1 < p < \infty$.

**Lemma 3.1.** *The transition densities with respect to the Lebesgue measure of the process $X^{\varepsilon,n}$ satisfy* (14).

**Proof.** Let us define a new $(d+1)$-dimensional process

$$d\widetilde{X}_t^{\varepsilon,n} = (0, b - n^{-1}DV(\widetilde{X}_t^{\varepsilon,n}, \omega))\,dt + (0, \sigma(\widetilde{X}_t^{\varepsilon,n}, \omega)\,dB_t) + n^{-1/2}dB_t'. \tag{16}$$

We point out that its generator on $C^2(\mathbb{R}^{d+1})$ can be rewritten in divergence form as

$$\widetilde{L}^{\varepsilon,n} = \left(\frac{e^{2V}}{2}\right)\sum_{i,j=0}^d \partial_i(e^{-2V}[a + H + n^{-1}I]_{ij}\partial_j)$$

(with the convention $a_{i0} = a_{0i} = H_{i0} = H_{0i} = 0$). From [12], this process admits transition densities with respect to the Lebesgue measure satisfying

$$\frac{1}{Ct^{(d+1)/2}}\exp\left(-\frac{C|y-x|^2}{t}\right) \le \widetilde{p}_{\varepsilon,n}(t,x,y) \le \frac{C}{t^{(d+1)/2}}\exp\left(-\frac{|y-x|^2}{Ct}\right) \tag{17}$$

for some constant $C > 0$ that only depends on $n, K, d$.

Let us now suppose that all the involved coefficients in (13) or (16) are smooth. For a $C_b^2(\mathbb{R}^{d+1})$ function $f$ (of class $C^2$ with bounded derivatives up to order 2), it is well-known that the mappings $(t,x) \in \mathbb{R} \times \mathbb{R}^{d+1} \mapsto \mathbb{E}_x[f(\widetilde{X}_t^{\varepsilon,n})]$ and $(t,x) \in \mathbb{R} \times \mathbb{R}^{d+1} \mapsto \mathbb{E}_x[f(X_t^{\varepsilon,n})]$ are at least of class $C^{1,2}(\mathbb{R} \times \mathbb{R}^{d+1})$ (see [7]) and are respectively classical solutions of the PDEs $\partial_t u = \widetilde{L}^{\varepsilon,n}u$ and $\partial_t u = \widetilde{L}^{\varepsilon,n}u + \varepsilon^{2\beta-\alpha}\partial_0 u$ with initial conditions $u(0,\cdot) = f$. Then, it is readily seen that the difference between these functions $(t,x) \mapsto \mathbb{E}_x[f(X_t^{\varepsilon,n})] - \mathbb{E}_{x+\varepsilon^{2\beta-\alpha}x_0}[f(\widetilde{X}_t^{\varepsilon,n})]$ is a classical solution of the PDE $\partial_t u = \widetilde{L}^{\varepsilon,n}u + \varepsilon^{2\beta-\alpha}\partial_0 u$ with initial condition $u(0,\cdot) = 0$. Thus, from the comparison theorem (see [17], Theorems 2.4 and 3.1, for a probabilistic proof of this fact) the functions coincide on $\mathbb{R}^{d+1}$. With density arguments, we establish $\mathbb{E}_x[f(X_t^{\varepsilon,n})] = \mathbb{E}_{x+\varepsilon^{2\beta-\alpha}x_0}[f(\widetilde{X}_t^{\varepsilon,n})]$ for each continuous bounded function $f$. The lemma then follows from (17).

If the coefficients in (13) or (16) are not smooth, the situation easily comes down to the previous case with the help of a regularization procedure. Details are however left to the reader. □

**Lemma 3.2.** *For each fixed* $\omega \in \Omega$ *and* $(s,x) \in \mathbb{R} \times \mathbb{R}^d$, *the process* $t \in \mathbb{R}_+ \mapsto (\varepsilon^{-\alpha}(t+s), \varepsilon^{-\beta}x + X_t^{\varepsilon,\tau_{s/\varepsilon^\alpha,x/\varepsilon^\beta}\omega})$, *starting from* $X_0^{\varepsilon,\tau_{s/\varepsilon^\alpha,x/\varepsilon^\beta}\omega} = 0$, *and the process* $t \in \mathbb{R}_+ \mapsto (\varepsilon^{-\alpha}(t+s), \varepsilon^{-\beta}X_{t+s}^{\varepsilon,\omega})$, *starting from* $X_s^{\varepsilon,\omega} = x$, *have the same law. As a consequence, for each fixed* $t \in \mathbb{R}_+$ *and* $\mathbf{f} \in C(\Omega)$, *the function* $P_t(\mathbf{f})$ *belongs to the space* $C(\Omega)$, *and*

$$P_t(\mathbf{f})(\tau_{\varepsilon^{2\beta-\alpha}s,x}\omega) = \mathbb{E}_x[f(\varepsilon^{2\beta-\alpha}(t+s), \overline{X}_{t+s}^\varepsilon, \omega)]. \tag{18}$$

**Proof.** For the sake of clarity and without loss of generality, let us drop the parameter $\varepsilon$ by choosing $\varepsilon = 1$. Fix $(s,x) \in \mathbb{R} \times \mathbb{R}^d$. It suffices to check that the process $(s+t, x+X_t^{\tau_{s,x}\omega})$, $X_t^{\tau_{s,x}\omega}$ starting from $0 \in \mathbb{R}^d$ at time $t=0$, and the process $(t+s, X_{t+s}^\omega)$, $X_{t+s}^\omega$ starting from $x \in \mathbb{R}^d$ at time $t=0$ (that is $X_s^\omega = x$), satisfy the same Itô equation, and therefore have the same law by virtue of law uniqueness for solutions of Itô's equations. As a consequence, for $\mathbf{f} \in C(\Omega)$, the function

$$(s,x) \mapsto \mathbb{E}[\mathbf{f}(Y_t(\tau_{s,x}\omega))] = \mathbb{E}_x[f(t+s, X_{t+s}^\omega, \omega)]$$

is continuous from the continuity of $\mathbf{f}$ and the Lipschitz properties of the coefficients $\mathbf{b}$ and $\boldsymbol{\sigma}$. □

## 4. Resolvent equation

Let us now investigate, for each $\lambda > 0$, the resolvent equation ($\mathbf{u}_\lambda$ is the unknown)

$$\lambda \mathbf{u}_\lambda - (\mathbf{L} + \theta D_t)\mathbf{u}_\lambda = \mathbf{h} \tag{19}$$

for a function $\mathbf{h} \in L^2(\Omega)$ and a function $\theta = \theta(\lambda)$ of the parameter $\lambda$. Note that the operator $\mathbf{L} + \theta D_t$ has not been rigorously defined yet but a complete description of all the involved operators and the meaning of "solution of (19)" are given thereafter. Actually, because of both time and space degeneracies, the generator of the process $Y^\varepsilon$ does not possess enough regularity to work on in a quite general way. However, for a suitable right-hand side $\mathbf{h}$, the function $\mathbf{u}_\lambda$ inherits some regularity properties that allow us to carry through the study of (19). This argument will be the guiding line of this section.

The following study is carried out for a quite general strictly positive function $\theta$ of the parameter $\lambda$. However, this result will only be used in Section 7 with $\theta = \theta(\lambda) = \lambda^{1-\alpha/(2\beta)}$. Roughly speaking, this function measures the difference of speed rates between the time and space variations.

### 4.1. Setup

In [18], it is proved that the unbounded operators $\widetilde{\mathbf{S}}$ and $D_t$ on $L^2(\Omega, \pi)$ (see the definition in Assumption 2.4) have a common spectral representation. This is due to the time independence of the coefficients $\widetilde{\mathbf{a}}$ and $\widetilde{\mathbf{b}}$. More precisely, we can find a spectral resolution of the identity $E$ on the Borelian subsets of $\mathbb{R}_+ \times \mathbb{R}$ such that

$$-\widetilde{\mathbf{S}} = \int_{\mathbb{R}_+ \times \mathbb{R}} x E(\mathrm{d}x, \mathrm{d}y) \quad \text{and} \quad -D_t = \int_{\mathbb{R}_+ \times \mathbb{R}} \mathrm{i}y E(\mathrm{d}x, \mathrm{d}y).$$

For any $\boldsymbol{\varphi}, \boldsymbol{\psi} \in L^2(\Omega)$, denote by $E_{\boldsymbol{\varphi},\boldsymbol{\psi}}$ the measure $E_{\boldsymbol{\varphi},\boldsymbol{\psi}} = (E\boldsymbol{\varphi}, \boldsymbol{\psi})_2$. Let $\mathbf{S}$ be the Friedrich extension of the operator defined on $\mathcal{C}$ by $(1/2)\mathrm{e}^{2\mathbf{V}} \sum_{i,j} D_i(\mathrm{e}^{-2\mathbf{V}} \mathbf{a}_{ij} D_j)$. For any $\boldsymbol{\varphi}, \boldsymbol{\psi} \in \mathcal{C}$, define

$$\langle \boldsymbol{\varphi}, \boldsymbol{\psi} \rangle_1 = -(\boldsymbol{\varphi}, \widetilde{\mathbf{S}}\boldsymbol{\psi})_2 = \int_{\mathbb{R}_+ \times \mathbb{R}} x E_{\boldsymbol{\varphi},\boldsymbol{\psi}}(\mathrm{d}x, \mathrm{d}y),$$

$$\langle \boldsymbol{\varphi}, \boldsymbol{\psi} \rangle_{1,S} = -(\boldsymbol{\varphi}, \mathbf{S}\boldsymbol{\varphi})_2,$$

and $\|\boldsymbol{\varphi}\|_1 = \langle \boldsymbol{\varphi}, \boldsymbol{\varphi} \rangle_1^{1/2}$, $\|\boldsymbol{\varphi}\|_{1,S} = \langle \boldsymbol{\varphi}, \boldsymbol{\varphi} \rangle_{1,S}^{1/2}$ the corresponding seminorms, whose kernels both match the $L^2(\Omega)$-sub-space of functions that are invariant under space translations (see Assumption 2.4). By virtue of assumption (7), these seminorms are equivalent

$$m\|\boldsymbol{\varphi}\|_1^2 \le \|\boldsymbol{\varphi}\|_{1,S}^2 \le M\|\boldsymbol{\varphi}\|_1^2. \tag{20}$$

Let $\mathbb{F}$ (respectively $\mathbb{H}$) be the Hilbert space that is the closure of $\mathcal{C}$ in $L^2(\Omega)$ with respect to the inner product $\varepsilon$ (resp. $\kappa$) defined on $\mathcal{C}$ by

$$\varepsilon(\boldsymbol{\varphi}, \boldsymbol{\psi}) = (\boldsymbol{\varphi}, \boldsymbol{\psi})_2 + \langle \boldsymbol{\varphi}, \boldsymbol{\psi} \rangle_1 + (D_t\boldsymbol{\varphi}, D_t\boldsymbol{\psi})_2$$

$$(\text{resp. } \kappa(\boldsymbol{\varphi}, \boldsymbol{\psi}) = (\boldsymbol{\varphi}, \boldsymbol{\psi})_2 + \langle \boldsymbol{\varphi}, \boldsymbol{\psi} \rangle_1).$$

Define the space $\mathbb{D}$ as the closure in $(L^2(\Omega), |\cdot|_2)$ of the subspace $\mathcal{D} = \{(-\widetilde{S})^{1/2}\boldsymbol{\varphi}; \boldsymbol{\varphi} \in \mathcal{C}\}$. For any $\boldsymbol{\varphi} \in \mathcal{C}$, define $\Phi((-\widetilde{S})^{1/2}\boldsymbol{\varphi}) = \boldsymbol{\sigma}^* D_x \boldsymbol{\varphi} \in (L^2(\Omega))^d$. Thanks to the description of the kernel of the semi-norm $\|\cdot\|_1$, note that $\Phi$ is well defined on $\mathcal{D}$. Furthermore $|\Phi((-\widetilde{S})^{1/2}\boldsymbol{\varphi})|_2^2 = -2(\boldsymbol{\varphi}, \mathbf{S}\boldsymbol{\varphi})_2$. From (20), $\Phi$ can be extended to the whole space $\mathbb{D}$ and this extension is a linear isomorphism from $\mathbb{D}$ into a closed subset of $(L^2(\Omega))^d$. Hence, for each function $\mathbf{u} \in \mathbb{H}$, we define $\nabla^\sigma \mathbf{u} = \Phi((-\widetilde{S})^{1/2}\mathbf{u})$ and this stands, in a way, for the gradient of $\mathbf{u}$ along the direction $\sigma$.

For each $\mathbf{f} \in L^2(\Omega)$ satisfying $\int_{\mathbb{R}_+ \times \mathbb{R}} \frac{1}{x} E_{\mathbf{f},\mathbf{f}}(\mathrm{d}x, \mathrm{d}y) < \infty$, we define

$$\|\mathbf{f}\|_{-1}^2 = \int_{\mathbb{R}_+ \times \mathbb{R}} \frac{1}{x} E_{\mathbf{f},\mathbf{f}}(\mathrm{d}x, \mathrm{d}y). \tag{21}$$

We point out that $\|\mathbf{f}\|_{-1} < \infty$ if and only if (see [14], for instance, for further details) there exists $C \in \mathbb{R}$ such that for any $\boldsymbol{\varphi} \in \mathcal{C}$, $(\mathbf{f}, \boldsymbol{\varphi})_2 \le C\|\boldsymbol{\varphi}\|_1$. For such a function $\mathbf{f}$, $\|\mathbf{f}\|_{-1}$ also matches the smallest $C$ satisfying this

inequality. Remark that $\|\mathbf{f}\|_{-1} < \infty$ implies $\pi(\mathbf{f}) = 0$. Denote by $\mathbb{H}_{-1}$ the closure of $L^2(\Omega)$ in $\mathbb{H}^*$ (topological dual of $\mathbb{H}$) with respect to the norm $\|\cdot\|_{-1}$.

Let us now focus on the antisymmetric part $\mathbf{H}$. We have

$$|(u, \mathbf{H}v)| \le (u, |\mathbf{H}|u)^{1/2}(v, |\mathbf{H}|v)^{1/2} \le C_1^H(u, \widetilde{\mathbf{a}}u)^{1/2}(v, \widetilde{\mathbf{a}}v)^{1/2}. \tag{22}$$

The second inequality follows from (8) and the first one is a general fact of linear algebra. We deduce

$$\forall \boldsymbol{\varphi}, \boldsymbol{\psi} \in \mathcal{C}, \quad \left(\frac{1}{2}\right)(\mathbf{H}D_x\boldsymbol{\varphi}, D_x\boldsymbol{\psi})_2 \le C_1^H \|\boldsymbol{\psi}\|_1 \|\boldsymbol{\varphi}\|_1. \tag{23}$$

For any $\boldsymbol{\varphi}, \boldsymbol{\psi} \in \mathcal{C}$, let us define

$$\left(\frac{1}{2}\right)(\mathbf{H}D_x\boldsymbol{\varphi}, D_x\boldsymbol{\psi})_2 = \mathbf{T}_H((-\widetilde{S})^{1/2}\boldsymbol{\varphi}, (-\widetilde{S})^{1/2}\boldsymbol{\psi}). \tag{24}$$

Note that, if $\boldsymbol{\varphi}, \boldsymbol{\varphi}' \in \mathcal{C}$ and $(-\widetilde{S})^{1/2}\boldsymbol{\varphi} = (-\widetilde{S})^{1/2}\boldsymbol{\varphi}'$, then $\widetilde{S}(\boldsymbol{\varphi} - \boldsymbol{\varphi}') = 0$ and, from Assumption 2.4, $D_x\boldsymbol{\varphi} = D_x\boldsymbol{\varphi}'$, in such a way that $\mathbf{T}_H((-\widetilde{S})^{1/2}\boldsymbol{\varphi}, (-\widetilde{S})^{1/2}\boldsymbol{\psi}) = \mathbf{T}_H((-\widetilde{S})^{1/2}\boldsymbol{\varphi}', (-\widetilde{S})^{1/2}\boldsymbol{\psi})$. Thus $\mathbf{T}_H$ is a well-defined antisymmetric bilinear operator on $\mathcal{D} \times \mathcal{D}$. From (23), it extends to an antisymmetric continuous bilinear form $\mathbf{T}_H$ on $\mathbb{D} \times \mathbb{D}$. Likewise, with the help of Assumption 2.3, we define the continuous bilinear forms $\mathbf{T}_a$, $\partial_t\mathbf{T}_a$, $\partial_t\mathbf{T}_H$, $\Lambda_s\mathbf{T}_a$, $\Lambda_s\mathbf{T}_a$ on $\mathbb{D} \times \mathbb{D} \subset L^2(\Omega, \pi) \times L^2(\Omega, \pi)$ as follows: $\forall \boldsymbol{\varphi}, \boldsymbol{\psi} \in \mathcal{C}$,

$$\left(\frac{1}{2}\right)(\mathbf{a}D_x\boldsymbol{\varphi}, D_x\boldsymbol{\psi})_2 = \mathbf{T}_a((-\widetilde{S})^{1/2}\boldsymbol{\varphi}, (-\widetilde{S})^{1/2}\boldsymbol{\psi}),$$

$$\left(\frac{1}{2}\right)(D_t\mathbf{a}D_x\boldsymbol{\varphi}, D_x\boldsymbol{\psi})_2 = \partial_t\mathbf{T}_a((-\widetilde{S})^{1/2}\boldsymbol{\varphi}, (-\widetilde{S})^{1/2}\boldsymbol{\psi}),$$

$$\left(\frac{1}{2}\right)(D_t\mathbf{H}D_x\boldsymbol{\varphi}, D_x\boldsymbol{\psi})_2 = \partial_t\mathbf{T}_H((-\widetilde{S})^{1/2}\boldsymbol{\varphi}, (-\widetilde{S})^{1/2}\boldsymbol{\psi}),$$

$$\left(\frac{1}{2}\right)(\Lambda_s\mathbf{a}D_x\boldsymbol{\varphi}, D_x\boldsymbol{\psi})_2 = \Lambda_s\mathbf{T}_a((-\widetilde{S})^{1/2}\boldsymbol{\varphi}, (-\widetilde{S})^{1/2}\boldsymbol{\psi}),$$

$$\left(\frac{1}{2}\right)(\Lambda_s\mathbf{H}D_x\boldsymbol{\varphi}, D_x\boldsymbol{\psi})_2 = \Lambda_s\mathbf{T}_H((-\widetilde{S})^{1/2}\boldsymbol{\varphi}, (-\widetilde{S})^{1/2}\boldsymbol{\psi}),$$

where, for any $s \in \mathbb{R}^*$, $\Lambda_s$ denotes the $L^2$-continuous difference operator (remind of the definition of $T_{s,0}$ in Section 2):

$$\forall \mathbf{f} \in L^2(\Omega), \quad \Lambda_s(\mathbf{f}) = \frac{(T_{s,0}\mathbf{f} - \mathbf{f})}{s}. \tag{25}$$

From Assumption 2.3, the norms of the forms $\Lambda_s\mathbf{T}_a$ and $\Lambda_s\mathbf{T}_H$ are uniformly bounded with respect to $s \in \mathbb{R}^*$ and the forms are weakly convergent respectively towards $\partial_t\mathbf{T}_a$ and $\partial_t\mathbf{T}_H$ as $s \to 0$.

Now, denote by $\mathcal{H}$ the subspace of $\mathbb{H}_{-1}$ whose elements satisfy the condition:

$$\exists C > 0, \forall s > 0 \text{ and } \forall \boldsymbol{\varphi} \in \mathcal{C}, \quad \langle \mathbf{h}, \Lambda_s\boldsymbol{\varphi} \rangle_{-1,1} \le C \|\boldsymbol{\varphi}\|_1. \tag{26}$$

For any $\mathbf{h} \in \mathcal{H}$, the smallest $C$ that satisfies such a condition is denoted $\|\mathbf{h}\|_T$. Then $\mathcal{H}$ is closed for the norm $\|\cdot\|_{\mathcal{H}} = \|\cdot\|_{-1} + \|\cdot\|_T$.

Finally, let us now extend the operator $\mathbf{L}$ defined on $\mathcal{C}$ by (12). For any $\lambda > 0$, consider the continuous bilinear form $\mathcal{B}_\lambda$ on $\mathbb{H} \times \mathbb{H}$ that coincides on $\mathcal{C} \times \mathcal{C}$ with

$$\forall \boldsymbol{\varphi}, \boldsymbol{\psi} \in \mathcal{C}, \quad \mathcal{B}_\lambda(\boldsymbol{\varphi}, \boldsymbol{\psi}) = \lambda(\boldsymbol{\varphi}, \boldsymbol{\psi})_2 + [\mathbf{T}_a + \mathbf{T}_H]((-\widetilde{S})^{1/2}\boldsymbol{\varphi}, (-\widetilde{S})^{1/2}\boldsymbol{\psi}).$$

Thanks to Assumption 2.3 and the antisymmetry of $\mathbf{H}$, this form is clearly coercive. Thus it defines a strongly continuous resolvent operator and consequently, a unique generator $\mathbf{L}$ associated to this resolvent operator. More precisely, $\boldsymbol{\varphi} \in \mathbb{H}$ belongs to $\mathrm{Dom}(\mathbf{L})$ if and only if $\mathcal{B}_\lambda(\boldsymbol{\varphi}, \cdot)$ is $L^2$-continuous. In this case, there exists $\mathbf{f} \in L^2(\Omega)$ such that $\mathcal{B}_\lambda(\boldsymbol{\varphi}, \cdot) = (\mathbf{f}, \cdot)_2$ and $\mathbf{L}\boldsymbol{\varphi}$ is equal to $\mathbf{f} - \lambda\boldsymbol{\varphi}$. It can be proved that this definition is independent of $\lambda > 0$ (see [11], Chapter 1, Section 2, for further details). Let us additionally mention that the adjoint operator $\mathbf{L}^*$ of $\mathbf{L}$ in $L^2(\Omega, \pi)$ can also be described through $\mathcal{B}_\lambda$. Indeed, $\mathrm{Dom}(\mathbf{L}^*) = \{\boldsymbol{\varphi} \in \mathbb{H}; \mathcal{B}_\lambda(\cdot, \boldsymbol{\varphi})$ is $L^2(\Omega)$-continuous$\}$. If $\boldsymbol{\varphi} \in \mathrm{Dom}(\mathbf{L}^*)$, there exists $\mathbf{f} \in L^2(\Omega)$ such that $\mathcal{B}_\lambda(\cdot, \boldsymbol{\varphi}) = (\mathbf{f}, \cdot)_2$ and $\mathbf{L}^*\boldsymbol{\varphi}$ is equal to $\mathbf{f} - \lambda\boldsymbol{\varphi}$.

**Remark 4.1.** *For each function $\boldsymbol{\varphi} \in \mathcal{C} \subset \mathbb{H}$, the application $\mathbf{L}\boldsymbol{\varphi}$ can be viewed as a function of $\mathbb{H}_{-1}$. Indeed, $\forall \boldsymbol{\psi} \in \mathcal{C}$, $(\mathbf{L}\boldsymbol{\varphi}, \boldsymbol{\psi})_2 = -[\mathbf{T}_a + \mathbf{T}_H]((-\widetilde{\mathbf{S}})^{1/2}\boldsymbol{\varphi}, (-\widetilde{\mathbf{S}})^{1/2}\boldsymbol{\psi}) \leq [M + C_1^H]\|\boldsymbol{\varphi}\|_1\|\boldsymbol{\psi}\|_1$. Hence, the application $\boldsymbol{\varphi} \mapsto \mathbf{L}\boldsymbol{\varphi} \in \mathbb{H}_{-1}$ can be extended to the whole space $\mathbb{H}$ so that, for each function $\mathbf{u} \in \mathbb{H}$, we can define $\mathbf{L}\mathbf{u}$ as an element of $\mathbb{H}_{-1}$ even if $\mathbf{u} \notin \mathrm{Dom}(\mathbf{L})$. The same properties hold for $\boldsymbol{\varphi} \in \mathbb{H} \mapsto \widetilde{\mathbf{S}}\boldsymbol{\varphi} \in \mathbb{H}_{-1}$ and $\boldsymbol{\varphi} \in \mathbb{H} \mapsto \mathbf{S}\boldsymbol{\varphi} \in \mathbb{H}_{-1}$.*

### 4.2. Existence of a solution

Let us now investigate the solvability of (19). Let us consider fixed parameters $\delta, \theta \geq 0$ and $\lambda > 0$ and introduce the bilinear form $B_{\lambda,\delta}^\theta$ on $\mathbb{F} \times \mathbb{F}$ that coincides on $\mathcal{C} \times \mathcal{C}$ with

$$B_{\lambda,\delta}^\theta(\boldsymbol{\varphi}, \boldsymbol{\psi}) = \lambda(\boldsymbol{\varphi}, \boldsymbol{\psi})_2 + \left(\frac{1}{2}\right)([\mathbf{a} + \mathbf{H}]D_x\boldsymbol{\varphi}, D_x\boldsymbol{\psi})_2 - \theta(D_t\boldsymbol{\varphi}, \boldsymbol{\psi})_2 + \left(\frac{\delta}{2}\right)(D_t\boldsymbol{\varphi}, D_t\boldsymbol{\psi})_2. \tag{27}$$

In what follows, the parameter $\delta$ is omitted each time that it is equal to 0. So the form $B_{\lambda,0}^\theta$ is simply denoted by $B_\lambda^\theta$.

The main result of this section is Proposition 4.2, whose proof can be found in [18], Proposition 5.4, without modification (except the value of $\theta$, which does not play a part). Let us however sum up the strategy. As explained in [18], the difficulty lies in the time degeneracy of the Dirichlet form $B_\lambda^\theta$, which satisfies no sector condition (even weak). To face this degeneracy, we add a viscosity parameter $\delta > 0$ and consider the form $B_{\lambda,\delta}^\theta$, which satisfies a weak sector condition. Thus, with the help of the Lax–Milgram theorem, it defines a resolvent operator that permits to solve ($\mathbf{u}_{\lambda,\delta}$ is the unknown)

$$\lambda\mathbf{u}_{\lambda,\delta} - (\mathbf{L} + \theta D_t)\mathbf{u}_{\lambda,\delta} - \left(\frac{\delta}{2}\right)D_t^2\mathbf{u}_{\lambda,\delta} = \mathbf{h}.$$

It then remains to let the parameter $\delta$ tend to 0. For this purpose, we have to establish estimates for the family $(|D_t\mathbf{u}_{\lambda,\delta}|_2)_\delta$ as $\delta$ goes to 0. This can be done when the right-hand side $\mathbf{h}$ is regular enough with respect to the time variable by using Assumption 2.3. Thus we state

**Proposition 4.2.** *Suppose that $\mathbf{h} \in L^2(\Omega) \cap \mathrm{Dom}(D_t)$ and $\mathbf{d} \in \mathcal{H}$. Then, for any fixed $\theta \geq 0$ and $\lambda > 0$, there exists a unique solution $\mathbf{u}_\lambda \in \mathbb{F}$ of the equation $\lambda\mathbf{u}_\lambda - \mathbf{L}\mathbf{u}_\lambda - \theta D_t\mathbf{u}_\lambda = \mathbf{h} + \mathbf{d}$, in the sense that $\forall \boldsymbol{\varphi} \in \mathbb{F}$, $B_\lambda^\theta(\mathbf{u}_\lambda, \boldsymbol{\varphi}) = (\mathbf{h}, \boldsymbol{\varphi})_2 + \langle \mathbf{d}, \boldsymbol{\varphi}\rangle_{-1,1}$. Moreover, $D_t\mathbf{u}_\lambda \in \mathbb{H}$ and we are provided with the following estimates, which only involve the constants $m, C_2^a, C_2^H$ of Assumption 2.3 (in particular, they depend neither on $\lambda$ nor on $\theta$),*

$$\lambda|\mathbf{u}_\lambda|_2^2 + m\|\mathbf{u}_\lambda\|_1^2 \leq \frac{|\mathbf{h}|_2^2}{\lambda} + \frac{\|\mathbf{d}\|_{-1}^2}{m}, \tag{28a}$$

$$\lambda|D_t\mathbf{u}_\lambda|_2^2 + m\|D_t\mathbf{u}_\lambda\|_1^2 \leq \frac{|D_t\mathbf{h}|_2^2}{\lambda} + 2\frac{\|\mathbf{d}\|_T^2}{m} + 2(C_2^a + C_2^H)^2\left(\frac{|\mathbf{h}|^2}{\lambda} + \frac{\|\mathbf{d}\|_{-1}^2}{m}\right) \Big/ m^2. \tag{28b}$$

*In the case $\mathbf{d} \in L^2(\Omega)$, then $\mathbf{u}_\lambda \in \mathrm{Dom}(\mathbf{L})$.*

*Finally, $\mathbf{u}_\lambda$ is the strong limit in $\mathbb{H}$ as $\delta$ goes to $0$ of the sequence $(\mathbf{u}_{\lambda,\delta})_{\lambda,\delta}$, where $\mathbf{u}_{\lambda,\delta}$ is the unique solution of the equation: $\forall \boldsymbol{\varphi} \in \mathbb{F}$, $B^\theta_{\lambda,\delta}(\mathbf{u}_{\lambda,\delta}, \boldsymbol{\varphi}) = (\mathbf{h}, \boldsymbol{\varphi})_2 + \langle \mathbf{d}, \boldsymbol{\varphi} \rangle_{-1,1}$, and the family $(D_t \mathbf{u}_{\lambda,\delta})_\delta$ is bounded in $L^2(\Omega)$.*

### 4.3. Control of the solution

We now aim at determining the asymptotic behavior, as $\lambda$ goes to $0$, of the solution $\mathbf{u}^i_\lambda$ (in the sense of Proposition 4.2) of the equation

$$\lambda \mathbf{u}^i_\lambda - (\mathbf{L} + \theta D_t)\mathbf{u}^i_\lambda = \mathbf{b}_i. \tag{29}$$

Since $\lambda$ is no more fixed, $\theta = \theta(\lambda)$ is not fixed either. More precisely, we are interested in two specific behaviors of the function $\theta = \theta(\lambda)$. We will focus on both possibilities $\lim_{\lambda \to 0} \theta(\lambda) = 0$ and $\lim_{\lambda \to 0} \theta(\lambda) = +\infty$, which respectively correspond to a small and big time/space evolution ratio. In both cases, we will show that $\lambda |\mathbf{u}^i_\lambda|^2_2 \to 0$ and that $(\nabla^\sigma \mathbf{u}^i_\lambda)_\lambda$ converges in $(L^2(\Omega))^d$ as $\lambda$ goes to $0$.

**Proposition 4.3.** *Let $(\mathbf{b}_\lambda)_{\lambda>0}$ be a family of functions in $\mathbb{H}_{-1} \cap L^2(\Omega)$ which is strongly convergent in $\mathbb{H}_{-1}$ to $\mathbf{b}_0$ and which is bounded in $\mathcal{H}$ (see definition (26)). Suppose either $\lim_{\lambda \to 0} \theta(\lambda) = 0$ or $\lim_{\lambda \to 0} \theta(\lambda) = +\infty$. Then the solution $\mathbf{u}_\lambda \in \mathbb{F}$ of the equation $\lambda \mathbf{u}_\lambda - \mathbf{L} \mathbf{u}_\lambda - \theta D_t \mathbf{u}_\lambda = \mathbf{b}_\lambda$ (in the sense of Proposition 4.2) satisfies:*

- *there exists $\boldsymbol{\eta} \in \mathbb{D}$ such that $(-\widetilde{\mathbf{S}})^{1/2} \mathbf{u}_\lambda \to \boldsymbol{\eta}$ as $\lambda$ goes to $0$ in $\mathbb{D}$,*
- *$\lambda |\mathbf{u}_\lambda|^2_2 \to 0$ as $\lambda$ goes to $0$.*

*The limit $\boldsymbol{\eta} \in \mathbb{D}$ does not depend on the function $\theta$ but only on its limit as $\lambda$ goes to $0$.*

The first step of the proof consists in investigating the case when the operator $\mathbf{L}$ is replaced by $\widetilde{\mathbf{S}}$ in (29). This situation is more convenient because of the common spectral decomposition of $\widetilde{\mathbf{S}}$ and $D_t$.

**Proposition 4.4.** *Suppose either $\lim_{\lambda \to 0} \theta(\lambda) = 0$ or $\lim_{\lambda \to 0} \theta(\lambda) = +\infty$. Let $(\mathbf{b}_\lambda)_{\lambda>0}$ be a family of functions in $\mathbb{H}_{-1}$ that is strongly convergent to $\mathbf{b}_0$ in $\mathbb{H}_{-1}$. Let $(\mathbf{v}_\lambda)_{\lambda>0}$ be a family of functions in $\mathbb{F}$ that solves the equation (for any $\lambda > 0$) $\lambda \mathbf{v}_\lambda - \widetilde{\mathbf{S}} \mathbf{v}_\lambda - \theta D_t \mathbf{v}_\lambda = \mathbf{b}_\lambda$ in the following sense,*

$$\forall \boldsymbol{\varphi} \in \mathbb{F}, \quad \lambda(\mathbf{v}_\lambda, \boldsymbol{\varphi})_2 + \langle \mathbf{v}_\lambda, \boldsymbol{\varphi} \rangle_1 - \theta(D_t \mathbf{v}_\lambda, \boldsymbol{\varphi})_2 = (\mathbf{b}_\lambda, \boldsymbol{\varphi})_2. \tag{30}$$

*Then there exists $\boldsymbol{\eta} \in \mathbb{D}$ such that $\lambda |\mathbf{v}_\lambda|^2_2 \to 0$ and $|(-\widetilde{\mathbf{S}})^{1/2} \mathbf{v}_\lambda - \boldsymbol{\eta}|_2 \to 0$ as $\lambda$ goes to $0$. The limit $\boldsymbol{\eta} \in \mathbb{D}$ does not depend on the function $\theta$ but only on its limit as $\lambda$ goes to $0$.*

**Proof.** *In what follows, the parameter $\lambda$ of $\theta(\lambda)$ is sometimes omitted when there is no possible confusion, but keep in mind that $\theta$ does depend on $\lambda$. From Lemmas 5.10 and 5.11 in [18], we can assume that, for any $\lambda > 0$, $\mathbf{b}_\lambda \in L^2(\Omega) \cap \text{Dom}(D_t) \cap \mathbb{H}_{-1}$ and converges to $\mathbf{b}_0 \in \mathbb{H}_{-1}$. From Proposition 4.2, we can suppose that $\mathbf{v}_\lambda \in \text{Dom}(\widetilde{\mathbf{S}})$. Remind that $-\widetilde{\mathbf{S}} = \int_{\mathbb{R}_+ \times \mathbb{R}} x E(\mathrm{d}x, \mathrm{d}y)$ and $-D_t = \int_{\mathbb{R}_+ \times \mathbb{R}} iy E(\mathrm{d}x, \mathrm{d}y)$. Choosing $\boldsymbol{\varphi} = \mathbf{v}_\lambda$ in (30), we have*

$$\lambda |\mathbf{v}_\lambda|^2_2 + \|\mathbf{v}_\lambda\|^2_1 = (\mathbf{b}_\lambda, \mathbf{v}_\lambda)_2 \leq C \|\mathbf{v}_\lambda\|_1 \leq C^2, \tag{31}$$

*where $C = \sup_{\lambda>0} \|\mathbf{b}_\lambda\|_{-1}$. Thus we can find $\boldsymbol{\eta} \in \mathbb{D}$ and a subsequence, still denoted by $(\mathbf{v}_\lambda)_\lambda$, such that $((-\widetilde{\mathbf{S}})^{1/2} \mathbf{v}_\lambda)_\lambda$ converges weakly in $L^2(\Omega)$ to $\boldsymbol{\eta}$.*

*Now we claim $\sup_{\lambda>0} \|\lambda \mathbf{v}_\lambda\|_{-1} < \infty$ and $\sup_{\lambda>0} \|\theta D_t \mathbf{v}_\lambda\|_{-1} < \infty$.*

$$\left| (\lambda \mathbf{v}_\lambda, \boldsymbol{\varphi})_2 \right| = \left| \int_{\mathbb{R}_+ \times \mathbb{R}} \lambda (\lambda + x + i\theta y)^{-1} \, \mathrm{d}E_{\mathbf{b}_\lambda, \boldsymbol{\varphi}}(\mathrm{d}x, \mathrm{d}y) \right|$$

$$\leq \left( \int_{\mathbb{R}_+ \times \mathbb{R}} \frac{\lambda^2}{x[(\lambda+x)^2 + \theta^2 y^2]} \, \mathrm{d}E_{\mathbf{b}_\lambda, \mathbf{b}_\lambda}(\mathrm{d}x, \mathrm{d}y) \right)^{1/2} \left( \int_{\mathbb{R}_+ \times \mathbb{R}} x \, \mathrm{d}E_{\boldsymbol{\varphi}, \boldsymbol{\varphi}}(\mathrm{d}x, \mathrm{d}y) \right)^{1/2}$$

$$\leq \sup_{\lambda>0}\left(\int_{\mathbb{R}_+\times\mathbb{R}} x^{-1}\,dE_{\mathbf{b}_\lambda,\mathbf{b}_\lambda}(dx,dy)\right)^{1/2}\|\boldsymbol{\varphi}\|_1$$

$$=\sup_{\lambda>0}\|\mathbf{b}_\lambda\|_{-1}\|\boldsymbol{\varphi}\|_1.$$

Since $\theta D_t\mathbf{v}_\lambda = \lambda\mathbf{v}_\lambda - \widetilde{\mathbf{S}}\mathbf{v}_\lambda - \mathbf{b}_\lambda$ and $\|\widetilde{\mathbf{S}}\mathbf{v}_\lambda\|_{-1}\leq\|\mathbf{v}_\lambda\|_1$ (cf. Remark 4.1), then $D_t\mathbf{v}_\lambda\in\mathbb{H}_{-1}$ and $\sup_{\lambda>0}\|\theta D_t\mathbf{v}_\lambda\|_{-1}<\infty$. So there exists a bounded family $(\mathbf{F}_\lambda)_{\lambda>0 0}$ of continuous linear forms on $\mathbb{D}\subset L^2(\Omega)$ such that $\forall\lambda>0$, $\forall\boldsymbol{\varphi}\in\mathcal{C}$, $\mathbf{F}_\lambda((-\widetilde{\mathbf{S}})^{1/2}\boldsymbol{\varphi})=\theta(\lambda)(D_t\mathbf{v}_\lambda,\boldsymbol{\varphi})_2$. Moreover, from (31), $(\lambda\mathbf{v}_\lambda)_\lambda$ converges to 0 in $L^2(\Omega)$ so that, $\forall\boldsymbol{\varphi}\in\mathcal{C}$

$$\mathbf{F}_\lambda((-\widetilde{\mathbf{S}})^{1/2}\boldsymbol{\varphi}) = (\lambda\mathbf{v}_\lambda,\boldsymbol{\varphi})_2+((-\widetilde{\mathbf{S}})^{1/2}\mathbf{v}_\lambda,(-\widetilde{\mathbf{S}})^{1/2}\boldsymbol{\varphi})_2-\langle\mathbf{b}_\lambda,\boldsymbol{\varphi}\rangle_{-1,1}$$
$$\to(\boldsymbol{\eta},(-\widetilde{\mathbf{S}})^{1/2}\boldsymbol{\varphi})_2-\langle\mathbf{b}_0,\boldsymbol{\varphi}\rangle_{-1,1} \tag{32}$$

as $\lambda$ goes to 0. So $(\mathbf{F}_\lambda)_{\lambda>0}$ is weakly convergent in $\mathbb{D}^*$ (topological dual of $\mathbb{D}$) to a limit denoted by $\mathbf{F}_0$ as $\lambda$ goes to 0.

We now aim at proving $\mathbf{F}_0(\boldsymbol{\eta})=0$. For $\lambda,\mu>0$, using the antisymmetry of the operator $D_t$ we obtain

$$\mathbf{F}_\lambda((-\widetilde{\mathbf{S}})^{1/2}\mathbf{v}_\mu)=\theta(\lambda)(D_t\mathbf{v}_\lambda,\mathbf{v}_\mu)_2=-\theta(\lambda)(D_t\mathbf{v}_\mu,\mathbf{v}_\lambda)_2=-\frac{\theta(\lambda)}{\theta(\mu)}\mathbf{F}_\mu((-\widetilde{\mathbf{S}})^{1/2}\mathbf{v}_\lambda). \tag{33}$$

If $\lim_{\lambda\to0}\theta(\lambda)=0$, we pass to the limit as $\lambda$ goes to 0 and we deduce $\mathbf{F}_0((-\widetilde{\mathbf{S}})^{1/2}\mathbf{v}_\mu)=0$. It just remains to pass to the limit as $\mu$ goes to 0 and it comes out $\mathbf{F}_0(\boldsymbol{\eta})=0$. Otherwise, if $\lim_{\lambda\to0}\theta(\lambda)=+\infty$, we first pass to the limit as $\mu$ goes to 0, then as $\lambda$ goes to 0, and we also obtain $\mathbf{F}_0(\boldsymbol{\eta})=0$.

Let us now establish the limit equation, which connects $\mathbf{F}_0$, $\boldsymbol{\eta}$ and $\mathbf{b}_0$. From (32), for any $\boldsymbol{\varphi}\in\mathcal{C}$

$$(\boldsymbol{\eta},(-\widetilde{\mathbf{S}})^{1/2}\boldsymbol{\varphi})_2-\mathbf{F}_0((-\widetilde{\mathbf{S}})^{1/2}\boldsymbol{\varphi})=\langle\mathbf{b}_0,\boldsymbol{\varphi}\rangle_{-1,1}. \tag{34}$$

This limit equation permits us to investigate uniqueness of the weak limit and, as a bypass, the last statement of Proposition 4.4. Indeed, let $\theta$ and $\theta'$ be two functions having the same limit, say $\lim_{\lambda\to0}\theta(\lambda)=\lim_{\lambda\to0}\theta'(\lambda)=0$ (the case $\lim_{\lambda\to0}\theta(\lambda)=\lim_{\lambda\to0}\theta'(\lambda)=+\infty$ is quite similar). Consider the two families of associated solutions $(\mathbf{v}_\lambda)_\lambda$ and $(\mathbf{v}'_\lambda)_\lambda$ of Eq. (30), and two possible weak limits $\mathbf{h}$ and $\boldsymbol{\eta}'$ of respectively two subsequences of $(\mathbf{v}_\lambda)_\lambda$ and $(\mathbf{v}'_\lambda)_\lambda$. Define the corresponding linear forms $(\mathbf{F}_\lambda)_\lambda$ and $(\mathbf{F}'_\lambda)_\lambda$ as well as their limits $\mathbf{F}_0$ and $\mathbf{F}'_0$ as described above. Then (34) provides us with the following equality:

$$\forall(-\widetilde{\mathbf{S}})^{1/2}\boldsymbol{\varphi}\in\mathbb{D},\quad(\boldsymbol{\eta}-\boldsymbol{\eta}',(-\widetilde{\mathbf{S}})^{1/2}\boldsymbol{\varphi})=[\mathbf{F}_0-\mathbf{F}'_0]((-\widetilde{\mathbf{S}})^{1/2}\boldsymbol{\varphi}). \tag{35}$$

As we proceeded for (33), we can establish

$$\mathbf{F}_\lambda((-\widetilde{\mathbf{S}})^{1/2}\mathbf{v}'_\mu)=-\frac{\theta(\lambda)}{\theta'(\mu)}\mathbf{F}'_\mu((-\widetilde{\mathbf{S}})^{1/2}\mathbf{v}_\lambda).$$

Now we pass to the limit as $\lambda$ goes to 0 along the first subsequence and then as $\mu$ goes to 0 along the second subsequence. We obtain $\mathbf{F}_0(\boldsymbol{\eta}')=0$. Reversing the roles of $\lambda$ and $\mu$, we also establish $\mathbf{F}'_0(\boldsymbol{\eta})=0$. Choose now $(-\widetilde{\mathbf{S}})^{1/2}\boldsymbol{\varphi}=\boldsymbol{\eta}-\boldsymbol{\eta}'$ in (35), this yields

$$|\boldsymbol{\eta}-\boldsymbol{\eta}'|_2^2=[\mathbf{F}_0-\mathbf{F}'_0](\boldsymbol{\eta}-\boldsymbol{\eta}')=0.$$

In particular, the weak convergence holds for the whole family $((-\widetilde{\mathbf{S}})^{1/2}\mathbf{v}_\lambda)_\lambda$ towards $\boldsymbol{\eta}$. Let us now tackle the strong convergence of $(\mathbf{v}_\lambda)_\lambda$. Choosing $\boldsymbol{\varphi}=\mathbf{v}_\lambda$ in (34), passing to the limit a $\lambda$ goes to 0 and using $\mathbf{F}_0(\mathbf{h})=0$, this yields

$$(\mathbf{h},\mathbf{h})_2=\lim_{\lambda\to0}\langle\mathbf{b}_0,\mathbf{v}_\lambda\rangle_{-1,1}=\lim_{\lambda\to0}\langle\mathbf{b}_\lambda,\mathbf{v}_\lambda\rangle_{-1,1}\overset{(30)}{=}\lim_{\lambda\to0}[\lambda|\mathbf{v}_\lambda|_2^2+\|\mathbf{v}_\lambda\|_1^2]. \tag{36}$$

In particular, $|\mathbf{h}|_2 = \lim_{\lambda \to 0} |(-\widetilde{\mathbf{S}})^{1/2}\mathbf{v}_\lambda|_2$. Thus, the convergence of the norms implies the strong convergence of the sequence $((-\widetilde{\mathbf{S}})^{1/2}\mathbf{v}_\lambda)_\lambda$ to $\boldsymbol{\eta}$ in $L^2(\Omega)$. As a bypass, (36) also implies the convergence of $(\lambda|\mathbf{v}_\lambda|_2^2)_\lambda$ to 0. $\qquad \square$

**Outline of the proof of Proposition 4.3.** Let us now briefly explain how to deduce Proposition 4.3 from Proposition 4.4. At first, note that $\mathbf{u}_\lambda = (\lambda - \mathbf{L} - D_t)^{-1}(\mathbf{b}_\lambda)$. The plan of attack then consists in defining the operator $P_\lambda \colon \mathcal{H} \to \mathcal{H}$ by $P_\lambda(\mathbf{b}) = (\mathbf{L} - \widetilde{\mathbf{S}})(\lambda - \widetilde{\mathbf{S}} - D_t)^{-1}(\mathbf{b})$ and in noting that

$$(\lambda - \mathbf{L} - D_t)^{-1} = (\lambda - \widetilde{\mathbf{S}} - D_t)^{-1}[\mathrm{I} - P_\lambda]^{-1}.$$

If $\|P_\lambda\|_{\mathcal{H} \to \mathcal{H}} < 1$, then the operator $[\mathrm{I} - P_\lambda]$ is invertible and $[\mathrm{I} - P_\lambda]^{-1}(\mathbf{b}_\lambda)$ converges in $\mathbb{H}_{-1}$. Thus Proposition 4.4 applies. Actually, $\|P_\lambda\|_{\mathcal{H} \to \mathcal{H}} \not< 1$ but this norm is finite. An iteration procedure shows that the general situation boils down to the particular case $\|P_\lambda\|_{\mathcal{H} \to \mathcal{H}} < 1$. Rigorous details can be found in [18], Proposition 5.12. $\qquad \square$

## 5. Itô's formula

Since the generator of the process $Y^\varepsilon$ is not regular enough to work on, we introduce regular approximations through viscosity methods. This allows us to get rid of the time degeneracy. Let us consider a standard 1-dimensional Brownian motion $\{B_t'; t \geq 0\}$ that is independent of $\{B_t; t \geq 0\}$ in such a way that $\{(B_t', B_t); t \geq 0\}$ is a standard $(d+1)$-dimensional Brownian motion. Let us then define the $(d+1)$-dimensional diffusion process $X^{\varepsilon,\delta}$, starting from 0, as the solution of the SDE (we still use the convention that $(x, X)$ stand for a $(d+1)$-dimensional vector where $x \in \mathbb{R}$ and $X \in \mathbb{R}^d$)

$$\overline{X}_t^{\varepsilon,\delta} = \int_0^t (\varepsilon^{2\beta-\alpha}, b(\overline{X}_r^{\varepsilon,\delta}, \omega)) \, \mathrm{d}r + \left(0, \int_0^t \sigma(\overline{X}_r^{\varepsilon,\delta}, \omega) \, \mathrm{d}B_r\right) + (\sqrt{\delta}B_t', 0). \tag{37}$$

The associated diffusion in random medium $Y^{\varepsilon,\delta}$ defined by $Y_t^{\varepsilon,\delta}(\omega) = \tau_{\overline{X}_t^{\varepsilon,\delta}}\omega$ is a $\Omega$-valued Markov process, which admits $\pi$ as invariant measure (similar to Section 3). So it defines a contraction semigroup on $L^2(\Omega, \pi)$. The associated (non-symmetric) Dirichlet form is given by (27) (with $\theta = \varepsilon^{2\beta-\alpha}$) with domain $\mathbb{F} \times \mathbb{F}$ and satisfies a weak sector condition (see [11], Chapter 1, Section 2, for the definition). The generator $\mathbf{L}^{\varepsilon,\delta}$ is defined on $\mathrm{Dom}(\mathbf{L}^{\varepsilon,\delta}) = \{\mathbf{u} \in \mathbb{F}; B_{\lambda,\delta}^\theta(\mathbf{u}, \cdot) \text{ is } L^2(\Omega)\text{-continuous}\}$ (see [11], Chapter 1, Section 2, for further details). It coincides on $\mathcal{C}$ with $\mathbf{L} + \theta D_t + (\delta/2)D_t^2$. Since $b$ and $\sigma$ are globally Lipschitz (Assumption 2.2), classical tools of SDE theory ensure that

$$\int_\Omega \mathbb{E}_0 \left[\sup_{0 \leq t \leq T} |(\varepsilon^{2\beta-\alpha}t, \overline{X}_t^\varepsilon) - X_t^{\varepsilon,\delta}|^2\right] \mathrm{d}\pi \to 0 \quad \text{as } \delta \text{ goes to } 0, \tag{38}$$

where both diffusions start from 0. The additional time regularity of this setting enables us to apply the Itô formula for suitable functions in the domain of $\mathbf{L}^{\varepsilon,\delta}$, typically functions $\mathbf{u}_{\lambda,\delta}$ given by Proposition 4.2. It remains to make the parameter $\delta$ vanish. Technical details are explained in [18].

**Theorem 5.1.** *Let* $\mathbf{f} \in L^2(\Omega)$ *and a family* $(\mathbf{u}_\lambda)_{\lambda>0}$ *in* $\mathbb{F}$ *such that:*

(1) $\forall \boldsymbol{\varphi} \in \mathbb{F}, \ B_\lambda^\theta(\mathbf{u}_\lambda, \boldsymbol{\varphi}) = (\mathbf{f}, \boldsymbol{\varphi})_2,$

(2) *for each* $\lambda > 0$, *there exists a sequence* $(\mathbf{u}_{\lambda,\delta})_{\delta>0}$ *in* $\mathbb{F}$ *that converges in* $\mathbb{H}$ *towards* $\mathbf{u}_\lambda$. *Moreover* $(\mathbf{u}_{\lambda,\delta})_{\delta>0}$ *satisfies* $B_{\lambda,\delta}^\theta(\mathbf{u}_{\lambda,\delta}, \cdot) = (\mathbf{f}, \cdot)_2,$

(3) *for each fixed* $\lambda > 0$, $(D_t\mathbf{u}_{\lambda,\delta})_\delta$ *is bounded in* $L^2(\Omega)$,

(4) *each function* $\mathbf{u}_{\lambda,\delta}$ *has continuous trajectories, that is, for* $\mu$ *almost every* $\omega \in \Omega$, *the function* $(t, x) \in \mathbb{R}^{d+1} \mapsto \mathbf{u}_{\lambda,\delta}(\tau_{t,x}\omega)$ *is continuous.*

*Then, choosing $\theta(\lambda) = \lambda^{1-\alpha/(2\beta)}$ and $\lambda = \varepsilon^{2\beta}$, $\mathbb{P}_\pi$ a.s., the following formula holds*

$$\mathbf{u}_{\varepsilon^{2\beta}}(Y_t^\varepsilon) = \mathbf{u}_{\varepsilon^{2\beta}}(Y_0^\varepsilon) + \int_0^t (\varepsilon^{2\beta}\mathbf{u}_{\varepsilon^{2\beta}} - \mathbf{f})(Y_r^\varepsilon)\,\mathrm{d}r + \int_0^t \nabla^\sigma \mathbf{u}_{\varepsilon^{2\beta}}^*(Y_r^\varepsilon)\,\mathrm{d}B_r,$$

*where $\mathbb{P}_\pi$ is the law of the process $Y^\varepsilon$ starting with initial distribution $\pi$ on $\Omega$.*

## 6. Ergodic theorem

We now focus on the ergodic properties of the family of processes $(Y^\varepsilon)_\varepsilon$. Basically, the difficulty lies here in the difference of evolution rates of the time and space variables. The strategy will consist in establishing an orthogonal decomposition of $L^2(\Omega)$ that allows, in a way, to separate the variables and to exploit separately their ergodic properties. Concerning the space variables, this is carried through from Assumption 2.4.

**Theorem 6.1.** *Let $\mathbf{f} \in L^1(\Omega)$. Then*

$$\mathbb{E}_\pi\left|\varepsilon^{2\beta}\int_0^{t/\varepsilon^{2\beta}} \mathbf{f}(Y_r^\varepsilon)\,\mathrm{d}r - t\pi(\mathbf{f})\right| \to 0 \quad \text{as } \varepsilon \text{ goes to } 0.$$

**Proof.** *Decomposition of the space $L^2(\Omega)$.* Let us first establish the following orthogonal decomposition of the space $L^2(\Omega)$

$$L^2(\Omega) = \mathbf{Inv} \oplus \mathrm{Closure}(\mathbf{R}),$$

where $\mathbf{R} = \{\widetilde{\mathbf{S}}\boldsymbol{\varphi}; \boldsymbol{\varphi} \in \mathrm{Dom}(\widetilde{\mathbf{S}})\} \cap \mathrm{Dom}(D_t)$ and $\mathbf{Inv} = \{\mathbf{f} \in L^2(\Omega); \forall x \in \mathbb{R}^d, T_{0,x}\mathbf{f} = \mathbf{f}\}$. It is readily seen that $\mathbf{Inv}$ and $\mathbf{R}$ are orthogonal and that $L^2(\Omega) \supset \mathbf{Inv} \oplus \mathrm{Closure}(\mathbf{R})$.

Conversely, let $\mathbf{f} \in L^2(\Omega)$ such that $\forall \widetilde{\mathbf{S}}\boldsymbol{\varphi} \in \mathbf{R}$, $(\mathbf{f}, \widetilde{\mathbf{S}}\boldsymbol{\varphi})_2 = 0$. Consider $\mathbf{h} \in L^2(\Omega) \cap \mathrm{Dom}(D_t)$. Hence Proposition 4.2 provides us with a family $(\mathbf{u}_\lambda)_\lambda \in \mathbb{F} \cap \mathrm{Dom}(\widetilde{\mathbf{S}})$ such that $B_\lambda^0(\mathbf{u}_\lambda, \boldsymbol{\varphi}) = (\mathbf{h}, \boldsymbol{\varphi})_2$ for any $\boldsymbol{\varphi} \in \mathbb{F}$. We emphasize that the equality

$$\lambda\mathbf{u}_\lambda - \widetilde{\mathbf{S}}\mathbf{u}_\lambda = \mathbf{h} \tag{39}$$

also holds in the $L^2$ sense, so that $\widetilde{\mathbf{S}}\mathbf{u}_\lambda = \lambda\mathbf{u}_\lambda - \mathbf{h} \in \mathbf{R}$. From (28a), $(\lambda\mathbf{u}_\lambda)_\lambda$ is bounded in $L^2(\Omega)$. Consider a weak limit $\mathbf{h}^* \in L^2(\Omega)$. Let us prove that $\widetilde{\mathbf{S}}\mathbf{h}^* = 0$. Let $\boldsymbol{\varphi} \in \mathrm{Dom}(\widetilde{\mathbf{S}})$. Integrate (39) against $\boldsymbol{\varphi}$, multiply by $\lambda$, pass to the limit as $\lambda$ goes to 0 and obtain $(\mathbf{h}^*, \widetilde{\mathbf{S}}\boldsymbol{\varphi})_2 = 0$. Since $\widetilde{\mathbf{S}}$ is self-adjoint, $\mathbf{h}^* \in \mathrm{Dom}(\widetilde{\mathbf{S}})$ and $\widetilde{\mathbf{S}}\mathbf{h}^* = 0$. From Assumption 2.4, $\mathbf{h}^* \in \mathbf{Inv}$. Let us additionally mention that $\mathbf{h}^*$ is the orthogonal projection of $\mathbf{h}$ into $\mathbf{Inv}$. Indeed, integrate (39) against $\boldsymbol{\varphi} \in \mathbf{Inv}$, note that $(\widetilde{\mathbf{S}}\mathbf{u}_\lambda, \boldsymbol{\varphi})_2 = 0$, pass to the limit as $\lambda$ goes to 0 and obtain $(\mathbf{h} - \mathbf{h}^*, \boldsymbol{\varphi})_2 = 0$. In particular, $(T_{0,x}\mathbf{h})^* = \mathbf{h}^*$ for $x \in \mathbb{R}^d$. Returning to our first objective, observe that $(\mathbf{f}, \mathbf{h})_2 = (\mathbf{f}, \lambda\mathbf{u}_\lambda - \widetilde{\mathbf{S}}\mathbf{u}_\lambda)_2 = (\mathbf{f}, \lambda\mathbf{u}_\lambda)_2$. Passing to the limit as $\lambda$ goes to 0, we obtain $(\mathbf{f}, \mathbf{h})_2 = (\mathbf{f}, \mathbf{h}^*)$ for each function $\mathbf{h} \in L^2(\Omega) \cap \mathrm{Dom}(D_t)$. We are now in position to conclude: for each $x \in \mathbb{R}^d$ and $\mathbf{h} \in L^2(\Omega) \cap \mathrm{Dom}(D_t)$,

$$(T_{0,x}\mathbf{f}, \mathbf{h})_2 = (\mathbf{f}, T_{0,-x}\mathbf{h})_2 = (\mathbf{f}, (T_{0,-x}\mathbf{h})^*)_2 = (\mathbf{f}, \mathbf{h}^*)_2 = (\mathbf{f}, \mathbf{h})_2.$$

Since $L^2(\Omega) \cap \mathrm{Dom}(D_t)$ is dense in $L^2(\Omega)$, we deduce $T_{0,x}\mathbf{f} = \mathbf{f}$. The decomposition follows.

*Particular case.* Let us first prove Theorem 6.1 in the case $\mathbf{f} = \mathbf{f}^* + \widetilde{\mathbf{S}}\boldsymbol{\varphi}$, where $\boldsymbol{\varphi} \in \mathrm{Dom}(\widetilde{\mathbf{S}})$, $\widetilde{\mathbf{S}}\boldsymbol{\varphi} \in \mathrm{Dom}(D_t)$ and $\mathbf{f}^* \in \mathbf{Inv}$. For each $\lambda > 0$, Proposition 4.2 provides us with a family $(\mathbf{f}_\lambda)_\lambda$ that satisfy (in the sense of Proposition 4.2)

$$\lambda\mathbf{f}_\lambda - \mathbf{L}\mathbf{f}_\lambda - \theta(\lambda)D_t\mathbf{f}_\lambda = \widetilde{\mathbf{S}}\boldsymbol{\varphi}. \tag{40}$$

Choose now $\lambda = \varepsilon^{2\beta}$ and $\theta(\lambda) = \lambda^{1-\alpha/(2\beta)}$. Theorem 5.1 yields

$$\varepsilon^{2\beta}(\mathbf{f}_{\varepsilon^{2\beta}}(Y^{\varepsilon}_{t/\varepsilon^{2\beta}}) - \mathbf{f}_{\varepsilon^{2\beta}}(\omega)) = \varepsilon^{2\beta} \int_0^{t/\varepsilon^{2\beta}} [\varepsilon^{2\beta} \mathbf{f}_{\varepsilon^{2\beta}} - \widetilde{\mathbf{S}}\boldsymbol{\varphi}](Y^{\varepsilon}_r) \, \mathrm{d}r + \varepsilon^{2\beta} \int_0^{t/\varepsilon^{2\beta}} \nabla^{\sigma} \mathbf{f}_{\varepsilon^{2\beta}}(Y^{\varepsilon}_r) \, \mathrm{d}B_r,$$

in such a way that, using the decomposition $\mathbf{f} = \mathbf{f}^{\ast} + \widetilde{\mathbf{S}}\boldsymbol{\varphi}$,

$$\varepsilon^{2\beta} \int_0^{t/\varepsilon^{2\beta}} \mathbf{f}(Y^{\varepsilon}_r) \, \mathrm{d}r = \varepsilon^{\alpha} \int_0^{t/\varepsilon^{\alpha}} \mathbf{f}^{\ast}(T_{r,0}\omega) \, \mathrm{d}r + \varepsilon^{2\beta} \int_0^{t/\varepsilon^{2\beta}} \varepsilon^{2\beta} \mathbf{f}_{\varepsilon^{2\beta}}(Y^{\varepsilon}_r) \, \mathrm{d}r$$

$$+ \varepsilon^{2\beta}(\mathbf{f}_{\varepsilon^{2\beta}}(Y^{\varepsilon}_{t/\varepsilon^{2\beta}}) - \mathbf{f}_{\varepsilon^{2\beta}}(\omega)) - \varepsilon^{2\beta} \int_0^{t/\varepsilon^{2\beta}} \nabla^{\sigma} \mathbf{f}_{\varepsilon^{2\beta}}(Y^{\varepsilon}_r) \, \mathrm{d}B_r.$$

Let us now prove that the above right-hand side behaves as its first term as $\varepsilon$ goes to 0. From Proposition 4.3 (note that $\widetilde{\mathbf{S}}\boldsymbol{\varphi} \in \mathbb{H}_{-1}$), we have $\lambda^2 |\mathbf{f}_{\lambda}|_2^2 + \lambda \|\mathbf{f}_{\lambda}\|_1^2 \to 0$ as $\lambda$ goes to 0. Thanks to the Jensen inequality, the invariance of the measure $\pi$ and Proposition 4.3 (note that $\widetilde{\mathbf{S}}\boldsymbol{\varphi} \in \mathbb{H}_{-1}$), we deduce

$$\mathbb{E}_{\pi}\left[\left(\varepsilon^{2\beta} \int_0^{t/\varepsilon^{2\beta}} \varepsilon^{2\beta} \mathbf{f}_{\varepsilon^{2\beta}}(Y^{\varepsilon}_r) \, \mathrm{d}r\right)^2\right] \leq t^2 \varepsilon^{4\beta} |\mathbf{f}_{\varepsilon^{2\beta}}|_2^2 \xrightarrow[\varepsilon \to 0]{} 0.$$

Similarly, $\mathbb{E}_{\pi}[(\varepsilon^{2\beta}(\mathbf{f}_{\varepsilon^{2\beta}}(Y^{\varepsilon}_{t/\varepsilon^{2\beta}}) - \mathbf{f}_{\varepsilon^{2\beta}}(\omega)))^2] \leq 2\varepsilon^{4\beta} |\mathbf{f}_{\varepsilon^{2\beta}}|_2^2$ and this latter quantity converges to 0 as $\varepsilon$ goes to 0. Concerning the martingale part, we just have to estimate its quadratic variations

$$\mathbb{E}_{\pi}\left[\left(\varepsilon^{2\beta} \int_0^{t/\varepsilon^{2\beta}} \nabla^{\sigma} \mathbf{f}_{\varepsilon^{2\beta}}(Y^{\varepsilon}_r) \, \mathrm{d}B_r\right)^2\right] = \mathbb{E}_{\pi}\left[\varepsilon^{4\beta} \int_0^{t/\varepsilon^{2\beta}} |\nabla^{\sigma} \mathbf{f}_{\varepsilon^{2\beta}}(Y^{\varepsilon}_r)|^2 \, \mathrm{d}r\right] = 2t\varepsilon^{2\beta} \|\mathbf{f}_{\varepsilon^{2\beta}}\|_1^2$$

to conclude that it also converges to 0. Thus it remains to study the convergence of the first term. The classical ergodic theory ensures that $\varepsilon^{\alpha} \int_0^{t/\varepsilon^{\alpha}} \mathbf{f}^{\ast}(T_{r,0}\omega) \, \mathrm{d}r$ converges in $L^2(\Omega, \pi)$ to a limit function that is invariant under time translations. Since $\mathbf{f}^{\ast}$ is already invariant under space translations, it is readily seen that this limit function is also invariant under space translations and consequently constant (see Definition 2.1). Thus it is equal to $\pi(\mathbf{f}^{\ast}) = \pi(\mathbf{f})$.

*General case.* Consider now a function $\mathbf{f} \in L^1(\Omega)$ and denote by $|\cdot|_1$ the $L^1(\Omega)$-norm. Obviously, we can find a sequence $(\mathbf{f}_n)_n \subset L^2(\Omega)$ that converges in $L^1(\Omega)$ towards $\mathbf{f}$. Then, using the invariance of the measure $\pi$,

$$\mathbb{E}_{\pi}\left|\varepsilon^{2\beta} \int_0^{t/\varepsilon^{2\beta}} \mathbf{f}(Y^{\varepsilon}_r) \, \mathrm{d}r - t\pi(\mathbf{f})\right|$$

$$\leq \mathbb{E}_{\pi}\left|\varepsilon^{2\beta} \int_0^{t/\varepsilon^{2\beta}} [\mathbf{f} - \mathbf{f}_n](Y^{\varepsilon}_r) \, \mathrm{d}r\right|$$

$$+ \mathbb{E}_{\pi}\left|\varepsilon^{2\beta} \int_0^{t/\varepsilon^{2\beta}} \mathbf{f}_n(Y^{\varepsilon}_r) \, \mathrm{d}r - t\pi(\mathbf{f}_n)\right| + t|\pi(\mathbf{f}) - \pi(\mathbf{f}_n)|$$

$$\leq 2t|\mathbf{f} - \mathbf{f}_n|_1 + \mathbb{E}_{\pi}\left|\varepsilon^{2\beta} \int_0^{t/\varepsilon^{2\beta}} \mathbf{f}_n(Y^{\varepsilon}_r) \, \mathrm{d}r - t\pi(\mathbf{f}_n)\right|.$$

It just remains to choose $n$ large enough to make the former term in the above right-hand side small and then to choose $\varepsilon$ small enough to treat the latter term. This completes the proof. $\qquad\square$

## 7. Invariance principle

*Notations.* From now on, fix $\theta(\lambda) = \lambda^{1-\alpha/(2\beta)}$. For $i \in \{1, \ldots, d\}$ and $\lambda > 0$, Proposition 4.2 provides us with a solution $\mathbf{u}_\lambda^i \in \mathbb{F} \cap \mathrm{Dom}(\mathbf{L})$ of the equation

$$\lambda \mathbf{u}_\lambda^i - \mathbf{L}\mathbf{u}_\lambda^i - \theta(\lambda) D_t \mathbf{u}_\lambda^i = \mathbf{b}_i,$$

and with a solution $\mathbf{c}_\lambda \in \mathbb{F} \cap \mathrm{Dom}(\mathbf{L})$ of the equation

$$\lambda \mathbf{c}_\lambda - \mathbf{L}\mathbf{c}_\lambda - \theta(\lambda) D_t \mathbf{c}_\lambda = \mathbf{c}.$$

An integration by parts proves that $\mathbf{b}, \mathbf{c}$ fulfill the assumptions of Proposition 4.3 (see [18], Lemma 5.14, for the detailed proof). Thus it makes sense to define $\boldsymbol{\xi}_i = \lim_{\lambda \to 0} \nabla^\sigma \mathbf{u}_\lambda^i$ and $\boldsymbol{\kappa} = \lim_{\lambda \to 0} \nabla^\sigma \mathbf{c}_\lambda$, where both limits are understood in the $L^2(\Omega)$ sense. Moreover, from Proposition 4.3, $\lambda |\mathbf{u}_\lambda|_2^2 + \lambda |\mathbf{c}_\lambda|_2^2 \to 0$ as $\lambda$ tends to 0.

Consider now a bounded continuous function $f \in C(\mathbb{R}^d)$ and the solution $z_{\varepsilon,\omega}$ of Eq. (2). The Feynman–Kac formula gives a probabilistic representation of $z_{\varepsilon,\omega}$ (cf. [17], Theorem 3.2 and Remark 3.3)

$$z_{\varepsilon,\omega}(t,x) = \mathbb{E}_x[f(X_t^\varepsilon)\exp(Q_t^\varepsilon)], \tag{41}$$

where $Q_t^\varepsilon = \int_0^t [\varepsilon^{-\beta}c + d](r/\varepsilon^\alpha, X_r^\varepsilon/\varepsilon^\beta, \omega)\,dr$. As guessed by the reader, the strategy consists in studying the convergence in law of the couple of processes $(X^\varepsilon, Q^\varepsilon)$ and then in establishing a uniform integrability argument.

Let us first tackle the convergence in law of the couple $(X^\varepsilon, Q^\varepsilon)$. Define $\overline{Q}_t^\varepsilon = \int_0^t [c + \varepsilon^\beta d](Y_r^\varepsilon)\,dr$ and remind, from Section 3, that $(X^\varepsilon, Q^\varepsilon)$ and $(\varepsilon^\beta \overline{X}_{./\varepsilon^{2\beta}}^\varepsilon, \varepsilon^\beta \overline{Q}_{./\varepsilon^{2\beta}}^\varepsilon)$ have the same law, where both processes $X^\varepsilon$ and $\overline{X}^\varepsilon$ start from 0.

Applying Theorem 5.1 to the functions $\mathbf{u}_{\varepsilon^{2\beta}}$ and $\mathbf{c}_{\varepsilon^{2\beta}}$ yields

$$\varepsilon^\beta \overline{X}_{t/\varepsilon^{2\beta}}^\varepsilon = H_t^\varepsilon + \varepsilon^\beta \int_0^{t/\varepsilon^{2\beta}} (\boldsymbol{\sigma} + \nabla^\sigma \mathbf{u}_{\varepsilon^{2\beta}}^*)(Y_r^\varepsilon)\,dB_r,$$

$$\varepsilon^\beta \overline{Q}_{t/\varepsilon^{2\beta}}^\varepsilon = G_t^\varepsilon + \varepsilon^{2\beta} \int_0^{t/\varepsilon^{2\beta}} \mathbf{d}(Y_r^\varepsilon)\,dr + \varepsilon^\beta \int_0^{t/\varepsilon^{2\beta}} \nabla^\sigma \mathbf{c}_{\varepsilon^{2\beta}}^*(Y_r^\varepsilon)\,dB_r,$$

where

$$H_t^\varepsilon = \varepsilon^{3\beta} \int_0^{t/\varepsilon^{2\beta}} \mathbf{u}_{\varepsilon^{2\beta}}(Y_r^\varepsilon)\,dr - \varepsilon^\beta \mathbf{u}_{\varepsilon^{2\beta}}(Y_{t/\varepsilon^{2\beta}}^\varepsilon) + \varepsilon^\beta \mathbf{u}_{\varepsilon^{2\beta}}(\omega),$$

$$G_t^\varepsilon = \varepsilon^{3\beta} \int_0^{t/\varepsilon^{2\beta}} \mathbf{c}_{\varepsilon^{2\beta}}(Y_r^\varepsilon)\,dr - \varepsilon^\beta \mathbf{c}_{\varepsilon^{2\beta}}(Y_{t/\varepsilon^{2\beta}}^\varepsilon) + \varepsilon^\beta \mathbf{c}_{\varepsilon^{2\beta}}(\omega).$$

For the sake of clarity, it is worth recalling that $\mathbb{P}_\pi$ denotes the law of the process $Y^\varepsilon$ starting with $\pi$ as initial distribution and that the measure $\pi$ is invariant for the process $Y^\varepsilon$. Let us now establish that the finite dimensional distributions of both processes $H^\varepsilon$ and $G^\varepsilon$ converge in $\mathbb{P}_\pi$-probability to 0. Using the Cauchy–Schwarz inequality and the invariance of the measure $\pi$, we obtain

$$\mathbb{E}_\pi[(H_t^\varepsilon)^2] \le 3(2 + t^2)\varepsilon^{2\beta}|u_{\varepsilon^{2\beta}}|_2^2,$$

and this latter quantity converges to 0 as $\varepsilon$ goes to 0. Likewise, $\mathbb{E}_\pi[(G_t^\varepsilon)^2]$ converges to 0 as $\varepsilon$ goes to 0. Then, from Theorem 6.1, the process $\varepsilon^{2\beta} \int_0^{\cdot/\varepsilon^{2\beta}} \mathbf{d}(Y_r^\varepsilon)\,dr$ converges, at least in probability, to the deterministic process $t \mapsto t\pi(\mathbf{d})$. Finally, we tackle the convergence of the martingale part of the process $(\varepsilon^\beta \overline{X}_{./\varepsilon^{2\beta}}^\varepsilon, \varepsilon^\beta \overline{Q}_{./\varepsilon^{2\beta}}^\varepsilon)$, which matches

$$\varepsilon^\beta \int_0^{t/\varepsilon^{2\beta}} [(\sigma + \nabla^\sigma \mathbf{u}_{\varepsilon^{2\beta}}^*), \nabla^\sigma \mathbf{c}_{\varepsilon^{2\beta}}^*](Y_r^\varepsilon)\,dB_r.$$

Clearly, under $\mathbb{P}_\pi$, the difference between this latter process and $\varepsilon^\beta \int_0^{t/\varepsilon^{2\beta}} [(\boldsymbol{\sigma} + \boldsymbol{\xi}^*), \boldsymbol{\kappa}^*](Y_r^\varepsilon) \, \mathrm{d}B_r$ vanishes as $\varepsilon$ goes to 0 in $\mathbb{P}_\pi$-quadratic mean. The quadratic variations are easily computed

$$\varepsilon^{2\beta} \int_0^{t/\varepsilon^{2\beta}} \begin{bmatrix} (\boldsymbol{\sigma} + \boldsymbol{\xi}^*)(\boldsymbol{\sigma} + \boldsymbol{\xi}^*)^* & (\boldsymbol{\sigma} + \boldsymbol{\xi}^*)\boldsymbol{\kappa} \\ \boldsymbol{\kappa}^*(\boldsymbol{\sigma} + \boldsymbol{\xi}^*)^* & \boldsymbol{\kappa}^*\boldsymbol{\kappa} \end{bmatrix} (Y_r^\varepsilon) \, \mathrm{d}r.$$

Using Theorem 6.1 again, these quadratic variations converge, at least in $\mathbb{P}_\pi$-probability, towards the deterministic process $t \mapsto At$, where the nonnegative symmetric matrix $A$ is given by

$$\overline{A} = \int_\Omega \begin{bmatrix} (\boldsymbol{\sigma} + \boldsymbol{\xi}^*)(\boldsymbol{\sigma} + \boldsymbol{\xi}^*)^*(\omega) & (\boldsymbol{\sigma} + \boldsymbol{\xi}^*)\boldsymbol{\kappa}(\omega) \\ \boldsymbol{\kappa}^*(\boldsymbol{\sigma} + \boldsymbol{\xi}^*)^*(\omega) & \boldsymbol{\kappa}^*\boldsymbol{\kappa}(\omega) \end{bmatrix} \, \mathrm{d}\pi(\omega).$$

To sum up, the finite dimensional distributions of the process $(\varepsilon^\beta \overline{X}_{t/\varepsilon^{2\beta}}^\varepsilon, \varepsilon^\beta \overline{Q}_{t/\varepsilon^{2\beta}}^\varepsilon)$ converge in law to the process $(0, D)t + \overline{A}^{1/2} B_t'$, where $D = \pi(\mathbf{d})$ and $B'$ is a $(d+1)$-dimensional Brownian motion. Actually, this convergence holds in the sense of weak convergence of processes in the space $C([0,T]; \mathbb{R}^{d+1})$ for each fixed $T > 0$. Section 8 is devoted to the proof of this fact.

## 8. Tightness

We keep the notations of Section 5. To prove the tightness of the underlying stochastic processes, we use Garsia, Rodemich and Rumsey inequality (GRR's inequality) (cf. [19]). More precisely, we follow [14], Section 3. Nevertheless, the required condition (43) is stronger than the one stated in [18], Section 9, that is $\mathbf{g} \in \mathbb{H}_{-1}$. But this methods provides us with a uniform integrability criterion needed in Section 9. Once again, we work on the viscosity approximations $Y^{\varepsilon,\delta}$ to avoid facing the time degeneracy of the generator of $Y^\varepsilon$.

From [14], Theorem 3.2, for a given function $\mathbf{g} \in L^\infty(\Omega)$, the function

$$\Gamma(\mathbf{g})(t, \omega) = \mathbb{E}\left[\exp\left(\int_0^t \mathbf{g}(Y_r^{\varepsilon,\delta}) \, \mathrm{d}r\right)\right]$$

belongs to $\mathrm{Dom}(\mathbf{L}^{\varepsilon,\delta})$ and satisfies $\partial_t \Gamma(\mathbf{g}) = \mathbf{L}^{\varepsilon,\delta}\Gamma(\mathbf{g}) + \mathbf{g}\Gamma(\mathbf{g})$ with initial condition $\Gamma(\mathbf{g})(0, \omega) = 1$. Then [14], Proposition 3.3, ensures that

$$|\Gamma(\mathbf{g})(t, \cdot)|_2^2 \leq \exp(2t \mathrm{Sp}(\mathbf{L}^{\varepsilon,\delta} + \mathbf{g})), \tag{42}$$

where $\mathrm{Sp}(\mathbf{L}^{\varepsilon,\delta} + \mathbf{g}) = \sup_{|\boldsymbol{\varphi}|_2=1}(\boldsymbol{\varphi}, [\mathbf{L}^{\varepsilon,\delta} + \mathbf{g}]\boldsymbol{\varphi})_2$ and the sup is taken over $\boldsymbol{\varphi} \in \mathrm{Dom}(\mathbf{L}^{\varepsilon,\delta}) \subset \mathbb{F}$. Thus, for any $\alpha > 0$, by the invariance of the measure $\pi$, we have

$$\mathbb{E}_\pi\left[\exp\left(\alpha\varepsilon^\beta \int_{s/\varepsilon^{2\beta}}^{t/\varepsilon^{2\beta}} \mathbf{g}(Y_r^{\varepsilon,\delta}) \, \mathrm{d}r\right)\right] = \int_\Omega \Gamma(\alpha\varepsilon^\beta \mathbf{g})(\varepsilon^{-2\beta}(t-s), \cdot) \, \mathrm{d}\pi$$

$$\leq |\Gamma(\alpha\varepsilon^\beta \mathbf{g})(\varepsilon^{-2\beta}(t-s), \cdot)|_2$$

$$\leq \exp((t-s)\mathrm{Sp}(\varepsilon^{-2\beta}\mathbf{L}^{\varepsilon,\delta} + \alpha\varepsilon^{-\beta}\mathbf{g})).$$

It just remains to give a bound for $\mathrm{Sp}(\varepsilon^{-2\beta}\mathbf{L}^{\varepsilon,\delta} + \alpha\varepsilon^{-\beta}\mathbf{g})$ and we want it to be independent of $\delta$ and $\varepsilon$. Suppose now that

$$\forall \boldsymbol{\varphi} \in \mathcal{C}, \quad |\boldsymbol{\varphi}|_2 = 1 \quad \Longrightarrow \quad (\mathbf{g}, \boldsymbol{\varphi}^2)_2 \leq C_{(43)} \|\boldsymbol{\varphi}\|_1 \tag{43}$$

for some constant $C_{(43)} > 0$. For any $\boldsymbol{\varphi} \in \mathbb{F}$ with $|\boldsymbol{\varphi}|_2 = 1$,

$$([\varepsilon^{-2\beta}\mathbf{L}^{\varepsilon,\delta} + \alpha\varepsilon^{-\beta}\mathbf{g}]\boldsymbol{\varphi}, \boldsymbol{\varphi})_2 \leq -\varepsilon^{-2\beta}m\|\boldsymbol{\varphi}\|_1^2 + \alpha\varepsilon^{-\beta}(\mathbf{g}, \boldsymbol{\varphi}^2)_2 - \left(\frac{\delta}{2}\right)|D_t \boldsymbol{\varphi}|_2^2$$

$$\leq -\varepsilon^{-2\beta} m \|\varphi\|_1^2 + \alpha \varepsilon^{-\beta} C_{(43)} \|\varphi\|_1$$

$$\leq \frac{\alpha^2 C_{(43)}^2}{4m}.$$

Let us additionally assume that, for each fixed $\omega$, the mapping $(t,x) \mapsto g(t,x,\omega)$ is globally Lipschitz. From (38), we can then pass to the limit as $\delta$ goes to 0 in the inequality

$$\mathbb{E}_\pi \left[ \exp \left( \alpha \varepsilon^\beta \int_{s/\varepsilon^{2\beta}}^{t/\varepsilon^{2\beta}} \mathbf{g}(Y_r^{\varepsilon,\delta}) \, dr \right) \right] \leq \exp \left( \frac{(t-s)\alpha^2 C_{(43)}^2}{4m} \right)$$

and get the bound

$$\mathbb{E}_\pi \left[ \exp \left( \alpha \varepsilon^\beta \int_{s/\varepsilon^{2\beta}}^{t/\varepsilon^{2\beta}} \mathbf{g}(Y_r^\varepsilon) \, dr \right) \right] \leq \exp \left( \frac{(t-s)\alpha^2 C_{(43)}^2}{4m} \right). \tag{44}$$

Let us now state GRR's inequality, whose proof can be found in [19] or [14], Proposition 3.1:

**Proposition 8.1 (Garsia–Rodemich–Rumsey's inequality).** *Let $p$ and $\Psi$ be strictly increasing continuous functions on $[0,+\infty[$ satisfying $p(0) = \Psi(0) = 0$ and $\lim_{t\to\infty} \Psi(t) = +\infty$. For given $T > 0$ and $f \in C([0,T];\mathbb{R}^d)$, suppose that there exists a finite $B$ such that;*

$$\int_0^T \int_0^T \Psi \left( \frac{|g(t) - g(s)|}{p(|t-s|)} \right) ds \, dt \leq B < \infty. \tag{45}$$

*Then, for all $0 \leq s \leq t \leq T$,*

$$|g(t) - g(s)| \leq 8 \int_0^{t-s} \Psi^{-1} \left( \frac{4B}{u^2} \right) dp(u). \tag{46}$$

Choose now $\Psi(t) = e^t - 1$, $\Psi^{-1}(t) = \ln(1+t)$ and $p(t) = \sqrt{t}$. As explained in [14], Section 3, (46) provides us with the following estimate of the continuity modulus:

$$\sup_{\substack{|t-s|\leq\delta \\ 0\leq s < t \leq T}} |g(t) - g(s)| \leq 8\sqrt{\delta} \ln(\delta^{-1}) \left[ \ln \left( 4 \int_0^T \int_0^T \exp \left( \frac{|g(t) - g(s)|}{\sqrt{|t-s|}} \right) ds \, dt \right) + 6 \right]. \tag{47}$$

Choosing $g(t) = \varepsilon^\beta \int_0^{t/\varepsilon^{2\beta}} \mathbf{g}(Y_r^\varepsilon) \, dr$, taking the expectation $\mathbb{E}_\pi$ and using Jensen's inequality, we obtain

$$\mathbb{E}_\pi \left[ \sup_{\substack{|t-s|\leq\delta \\ 0\leq s < t \leq T}} \left| \varepsilon^\beta \int_{s/\varepsilon^{2\beta}}^{t/\varepsilon^{2\beta}} \mathbf{g}(Y_r^\varepsilon) \, dr \right| \right]$$

$$\leq 8\sqrt{\delta} \ln(\delta^{-1}) \left[ \ln \left( 4 \int_0^T \int_0^T \mathbb{E}_\pi \left[ \exp \left( \frac{|\varepsilon^\beta \int_{s/\varepsilon^{2\beta}}^{t/\varepsilon^{2\beta}} \mathbf{g}(Y_r^\varepsilon) \, dr|}{\sqrt{|t-s|}} \right) \right] ds \, dt \right) + 6 \right].$$

Gathering (44) with the above estimate, we finally obtain the desired continuity modulus estimate

$$\mathbb{E}_\pi \left[ \sup_{\substack{|t-s|\leq\delta \\ 0\leq s < t \leq T}} \left| \varepsilon^{-\beta} \int_{s/\varepsilon^{2\beta}}^{t/\varepsilon^{2\beta}} \mathbf{g}(Y_r^\varepsilon) \, dr \right| \right] \leq C_T \sqrt{\delta} \ln(\delta^{-1}) \tag{48}$$

for some positive constant $C_T$. Note that both functions $\mathbf{b}$ and $\mathbf{c}$ satisfy (43). Indeed, for any $\boldsymbol{\varphi} \in \mathcal{C}$ (the case $\mathbf{g} = \mathbf{b}$ is similar)

$$
\begin{aligned}
(\mathbf{c}, \boldsymbol{\varphi}^2)_2 &= \sum_{i,j} (e^{2\mathbf{V}} D_i (e^{-2\mathbf{V}} \widetilde{\boldsymbol{\sigma}}_{ij} \mathbf{f}_j), \boldsymbol{\varphi}^2)_2 = -\sum_{i,j} (\widetilde{\boldsymbol{\sigma}}_{ij} \mathbf{f}_j, D_i(\boldsymbol{\varphi}^2))_2 \\
&= -2 \sum_{i,j} (\boldsymbol{\varphi} \mathbf{f}_j, \widetilde{\boldsymbol{\sigma}}_{ij} D_i \boldsymbol{\varphi})_2 \le 4 \|\mathbf{f}\|_{L^\infty(\Omega)} |\boldsymbol{\varphi}|_2 \|\boldsymbol{\varphi}\|_1.
\end{aligned}
$$

In particular, (48) holds for $\mathbf{g} = \mathbf{b}$ and $\mathbf{g} = \mathbf{c}$. Note that, from the boundedness of $\mathbf{d}$ and $\boldsymbol{\sigma}$, the Burkholder–Davis–Gundy and the Kolmogorov criterion, the tightness of the processes $\varepsilon^{2\beta} \int_0^{\cdot/\varepsilon^{2\beta}} \mathbf{d}(Y_r^\varepsilon) \, dr$ and $\varepsilon^\beta \int_0^{\cdot/\varepsilon^{2\beta}} \boldsymbol{\sigma}(Y_r^\varepsilon) \, dB_r$ raises no particular difficulty. As a consequence, the tightness of the processes $\varepsilon^\beta \overline{X}_{\cdot/\varepsilon^{2\beta}}^\varepsilon$ and $\varepsilon^\beta \overline{Q}_{\cdot/\varepsilon^{2\beta}}^\varepsilon$ is proved.

## 9. Girsanov's transform and limit equation

To determine the limit of $z_{\varepsilon,\omega}(t,x)$ (see (41)), we now aim at applying the convergence in law of the couple $(\varepsilon^\beta \overline{X}_{t/\varepsilon^{2\beta}}^\varepsilon, \varepsilon^\beta \overline{Q}_{t/\varepsilon^{2\beta}}^\varepsilon)$ to the function $g(x_1, \dots, x_{d+1}) = f(x_1, \dots, x_d) e^{x_{d+1}}$ for some continuous bounded function $f \in C_b(\mathbb{R}^d; \mathbb{R})$. This requires a uniform integrability argument, which is derived from (44), applied with $\mathbf{g} = \mathbf{c}$. Since $\mathbf{d}$ is bounded, the random variable $\exp(\varepsilon^{2\beta} \int_0^{t/\varepsilon^{2\beta}} \mathbf{d}(Y_r^\varepsilon) \, dr)$ is bounded too. So the random variable $f(\varepsilon^\beta \overline{X}_{t/\varepsilon^{2\beta}}^\varepsilon) \exp(\varepsilon^\beta \overline{Q}_{t/\varepsilon^{2\beta}}^\varepsilon)$ is uniformly integrable and thus converges in $\pi$-probability (cf. Section 7) to $\mathbb{E}_0[g((0,D)t + \overline{A}^{1/2} B_t')]$. We now aim at finding a PDE that characterizes this limit. Let us consider the $\mathbb{R}^d \times \mathbb{R}$-de-composition in blocs of $\overline{A}$ and $\overline{A}^{1/2}$

$$
\overline{A} = \left( \begin{array}{c|c} \overline{A}_{11} & \overline{A}_{12} \\ \hline \overline{A}_{21} & \overline{A}_{22} \end{array} \right), \qquad \overline{A}^{1/2} = \left( \begin{array}{c|c} (\overline{A}^{1/2})_{11} & (\overline{A}^{1/2})_{12} \\ \hline (\overline{A}^{1/2})_{21} & (\overline{A}^{1/2})_{22} \end{array} \right).
$$

By the Girsanov transform, we can define a new probability $\widetilde{\mathbb{P}}$ on $C([0,T]; \mathbb{R}^{d+1})$ as follows

$$
\forall t < T, \quad \left. \frac{d\widetilde{\mathbb{P}}}{d\mathbb{P}} \right|_{\mathcal{F}_t} = \exp\left( \left( \overline{A}_{21}^{1/2} \quad \overline{A}_{22}^{1/2} \right) B_t' - \overline{A}_{22} \frac{t}{2} \right),
$$

where $\mathcal{F}_t$ denotes the natural filtration of the Brownian motion. Under $\widetilde{\mathbb{P}}$, $\widetilde{B}_t = B_t' - (\overline{A}_{12}^{1/2}, \overline{A}_{22}^{1/2})^* t$ is a Brownian motion. Let us now rewrite $\mathbb{E}_0[g((0,D)t + \overline{A}^{1/2} B_t')]$ in terms of $\widetilde{B}$

$$
\begin{aligned}
&\mathbb{E}_0[g((0,D)t + \overline{A}^{1/2} B_t')] \\
&= \mathbb{E}[f((\overline{A}_{11}^{1/2}, \overline{A}_{12}^{1/2}) B_t') \exp(Dt + (\overline{A}_{21}^{1/2}, \overline{A}_{22}^{1/2}) B_t')] \\
&= \widetilde{\mathbb{E}}_0 \left[ f((\overline{A}_{11}^{1/2}, \overline{A}_{12}^{1/2}) \widetilde{B}_t + t(\overline{A}_{11}^{1/2} \overline{A}_{12}^{1/2} + \overline{A}_{12}^{1/2} \overline{A}_{22}^{1/2})) \exp\left( Dt + \overline{A}_{22} \frac{t}{2} \right) \right] \\
&= \widetilde{\mathbb{E}}[f((\overline{A}_{11}^{1/2} \quad \overline{A}_{12}^{1/2}) \widetilde{B}_t + Ct)] e^{Ut},
\end{aligned}
$$

where

$$
C = \int_\Omega (\boldsymbol{\sigma} + \boldsymbol{\xi}^*) \boldsymbol{\kappa} \, d\pi \quad \text{and} \quad U = \int_\Omega \left( \frac{\boldsymbol{\kappa}^* \boldsymbol{\kappa}}{2} + \mathbf{d} \right) d\pi. \tag{49}
$$

We point out that the quadratic variations of $(\overline{A}_{11}^{1/2}\overline{A}_{12}^{1/2})\widetilde{B}_t$ are equal to

$$A = \overline{A}_{11} = \int_\Omega (\boldsymbol{\sigma} + \boldsymbol{\xi}^*)(\boldsymbol{\sigma} + \boldsymbol{\xi}^*)^* \, d\pi. \tag{50}$$

Hence there exists a standard $d$-dimensional Brownian motion $\overline{B}_t$ such that

$$\mathbb{E}_0[g((0,D)t + \overline{A}^{1/2}B_t')] = \overline{\mathbb{E}}_0[f(Ct + A^{1/2}\overline{B}_t)e^{Ut}].$$

To sum up, we identified the limit in $\pi$ probability of

$$z_{\varepsilon,\omega}(0,t) = \mathbb{E}_0[f(X_t^\varepsilon)\exp(Q_t^\varepsilon)] = \mathbb{E}_0[f(\varepsilon^\beta \overline{X}_{t/\varepsilon^{2\beta}}^\varepsilon)\exp(\varepsilon^\beta \overline{Q}_{t/\varepsilon^{2\beta}}^\varepsilon)]$$

as $\overline{\mathbb{E}}_0[f(Ct + A^{1/2}\overline{B}_t)e^{Ut}]$.

Let us now determine this limit when the starting point is not 0 but $x \in \mathbb{R}^d$:

$$\mathbb{E}_x[f(X_t^{\varepsilon,\omega})\exp(Q_t^{\varepsilon,\omega})] \overset{\text{Lemma 3.2}}{=} \mathbb{E}_0[f(x + X_t^{\varepsilon,\tau_{(0,x/\varepsilon^\beta)}\omega})\exp(Q_t^{\varepsilon,\tau_{(0,x/\varepsilon^\beta)}\omega})]$$

$$\overset{\text{in law with respect to }\mu}{=} \mathbb{E}_0[f(x + X_t^{\varepsilon,\omega})\exp(Q_t^{\varepsilon,\omega})]$$

$$\underset{\varepsilon \to 0}{\overset{\pi \text{ prob}}{\longrightarrow}} \overline{\mathbb{E}}_0[f(x + Ct + A^{1/2}\overline{B}_t)e^{Ut}].$$

To complete the proof of Theorem 2.6, it just remains to explain why the coefficients $A, C, U$ only depend on the 3 cases $\alpha - 2\beta < 0, \alpha - 2\beta > 0$ or $\alpha - 2\beta = 0$. Since $\boldsymbol{\xi}$ and $\boldsymbol{\kappa}$ only depend on the case $\alpha - 2\beta < 0$, $\alpha - 2\beta > 0$ and $\alpha - 2\beta = 0$ (see Proposition 4.3 for the first two cases and [18] for the last one), the same property holds for $A$, $C$ and $U$, which can be expressed in terms of $\boldsymbol{\xi}$ and $\boldsymbol{\kappa}$ (see (49) and (50)). As a consequence, this completes the proof of Theorem 2.6.